\documentclass[11pt,reqno]{amsart}
\usepackage{ulem,soul}
\usepackage{amsmath,amsfonts,amssymb,mathrsfs,graphicx,subfigure}
\usepackage{color,psfrag,graphicx,subfigure}
\usepackage{algorithm,algpseudocode,enumerate}
\usepackage[subnum]{cases}
\usepackage{lineno,empheq}
\usepackage[margin=1.1in]{geometry}    

\usepackage{appendix}

 \newcommand{\nn}{\nonumber}

\newtheorem{theorem}{Theorem}[section]
\newtheorem{remark}[theorem]{Remark}

\newtheorem{corollary}[theorem]{Corollary}
\newtheorem{proposition}[theorem]{Proposition}
\newtheorem{definition}[theorem]{Definition}
\newcommand{\secref}[1]{\S\ref{#1}}

\usepackage[pdftex,colorlinks,bookmarksopen,bookmarksnumbered,citecolor=red,urlcolor=red]{hyperref}

\newcommand{\ld}{{\sigma}}
\newcommand{\ldsq}{{\sigma^2}}

\newcommand{\PL}{\textmd{I}}

 \newcommand{\EE}{\mathbb{E}}

 \newcommand{\Real}{\mathbf{R}}

\def\bn{\boldsymbol{n}}
\def\leq{\leqslant}
\def\geq{\geqslant}

\def\Lo{\mathcal{L}}

\newcommand{\wt}[1]{\widetilde{#1}}

\newcommand{\bh}{\bar{h}}

\newcommand{\eps}{\epsilon}

\newcommand{\norm}[1]{\lVert#1\rVert}

 \newcommand{\set}[1]{ \left\{#1\right\}}

\newcommand{\ud}{\,\mathrm{d}}
\newcommand{\rd}{\mathrm{d}}

\newcommand{\dom}{\mathcal{D}}

\newcommand{\Od}{\mathcal{O}}

\newcommand{\rex}{\mathrm{lay}}

\begin {document}

\title[]
{Estimation of exciton diffusion lengths of organic semiconductors in random domains
\vrule height
15pt width 0pt}

%
%
%

\date{\today}
\begin {abstract}
Exciton diffusion length plays a vital role in the function of opto-electronic devices. Oftentimes, the domain occupied by an organic semiconductor is subject to surface measurement error. In many experiments, photoluminescence over the domain is measured and used as the observation data to estimate this length parameter in an inverse manner based on the least square method. However, the result is sometimes found to be sensitive to the surface geometry of the domain. In this paper, we employ a random function representation for the uncertain surface of the domain. After non-dimensionalization, the forward model becomes a diffusion-type equation over the domain whose geometric boundary is subject to small random perturbations. We propose an asymptotic-based method as an approximate forward solver whose accuracy is justified both theoretically and numerically. It only requires solving several deterministic problems over a fixed domain. Therefore, for the same accuracy requirements we tested here, the running time of our approach is more than one order of magnitude smaller than that of directly solving the original stochastic boundary-value problem by the stochastic collocation method. In addition, from numerical results, we find that the correlation length of randomness is important to determine whether a 1D reduced model is a good surrogate.
\end {abstract}


\begin{center}
{\bf \Large Estimation of exciton diffusion lengths of organic semiconductors in random domains}

\

\

Jingrun Chen \footnote{email:
{jingrunchen@suda.edu.cn}.}
\par
\
{Mathematical Center for Interdisciplinary Research and
\\
School of Mathematical Sciences,
\\
Soochow University, Suzhou, China}
\par

\

\par

Ling Lin \footnote{email:
{linling27@mail.sysu.edu.cn}.}
\par
\
{School of Mathematics,
Sun Yat-Sen University,
\par
Guang Zhou, China}
\par

\
\par
Zhiwen Zhang
\footnote{email : zhangzw@hku.hk.  Corresponding author}

Department of Mathematics, The University of Hong Kong,
\par
 Pokfulam, Hong Kong SAR

\par
\
\par
Xiang Zhou \footnote{email: xiang.zhou@cityu.edu.hk.}\par
Department of Mathematics,
City University of Hong Kong\par
Tat Chee Ave, Kowloon,  Hong Kong SAR

\end{center}

\date{\today}
\section*{abstract}
Exciton diffusion length plays a vital role in the function of opto-electronic devices. Oftentimes, the domain occupied by an organic semiconductor is subject to surface measurement error. In many experiments, photoluminescence over the domain is measured and used as the observation data to estimate this length parameter in an inverse manner based on the least square method. However, the result is sometimes found to be sensitive to the surface geometry of the domain. In this paper, we employ a random function representation for the uncertain surface of the domain. After non-dimensionalization, the forward model becomes a diffusion-type equation over the domain whose geometric boundary is subject to small random perturbations. We propose an asymptotic-based method as an approximate forward solver whose accuracy is justified both theoretically and numerically. It only requires solving several deterministic problems over a fixed domain. Therefore, for the same accuracy requirements we tested here, the running time of our approach is more than one order of magnitude smaller than that of directly solving the original stochastic boundary-value problem by the stochastic collocation method. In addition, from numerical results, we find that the correlation length of randomness is important to determine whether a 1D reduced model is a good surrogate for the 2D model.

\medskip
{{\bf Subject class[2000]} {34E05, 35C20, 35R60, 58J37, 65C99}}
\par
{{\bf Keywords:} exciton diffusion, random domain, asymptotic methods, uncertainty qualification, organic semiconductor.}

\section{Introduction}
\noindent
From a practical perspective, measurement error or insufficient data in many problems
inevitably introduces uncertainty, which however has been overlooked for a long
time. In materials science, recent adventure in manufacturing has reduced the device dimension
from macroscropic/mesoscropic scales to nanoscale, in which the uncertainty becomes
important \cite{Bejan:2000}. In the field of organic opto-electronics, such as organic
light-emitting diodes (LEDs) and organic photovoltaics, a surge of interest has occurred
over the past few decades, due to major advancements in material design, which led to a
significant boost in the materials performance \cite{PopeSwenberg:1999, MyersXue:2012, SuLanWei:2012}. These materials are
carbon-based compounds with other elements like N, O, H, S, and P, and can be classified
into small molecules, oligomers, and polymers with atomic mass units ranging from several hundreds to
at least several thousands and conjugation length ranging from a few nanometers to hundreds of nanometers \cite{Forrest:2004, MyersXue:2012}.

At the electronic level, exciton, a bound electron-hole pair, is the elementary energy carrier,
which does not carry net electric charge. The characteristic distance that an exciton travels
during its lifetime is defined as the exciton diffusion length, which plays a critical  role in the
function of opto-electronical devices. A small diffusion length in organic photovoltaics limits
the dissociation of excitons into free charges \cite{TeraoSasabeAdachi:2007, MenkeLuhmanHolmes:2012},
while a large diffusion length in organic LEDs may limit luminous efficiency if excitons diffuse to non-radiative
quenching sites \cite{Antoniadis:1994}. Generally, there are two types of experimental methods to
measure exciton diffusion length:
photoluminescence quenching measurement, including steady-state
and time-resolved photoluminescence surface quenching, time-resolved photoluminescence
bulk quenching, and exciton-exciton annihilation \cite{Linetal:2013}, and photocurrent spectrum measurement \cite{PetterssonRomanInganas:1999}.
Exciton generation, diffusion, dissociation, recombination, exciton-exciton annihilation,
and exciton-environment interaction, are the typical underlying processes.
Accordingly, two types of models are used to describe exciton diffusion, either
differential equation based or stochastic process based.
The connections between these
models are systematically discussed in \cite{Chen:2016}.

We focus on the differential equation  model in this paper.
Accordingly, the device used in the experiment
includes two layers of organic materials. One layer of material is called donor and the other is
called acceptor or quencher due to the difference of their chemical properties. A typical bilayer
structure is illustrated in Figure \ref{fig:dom}. These materials are thin films with thicknesses
ranging from tens of nanometers to hundreds of nanometers along the $x$ direction and in-plane
dimensions up to the macroscopic scale. Under the illumination of solar lights, excitons are generated
in the donor layer, and then diffuse. Due to the exciton-environment interaction, some excitons
die out and emit photons which contribute to the photoluminescence. The donor-acceptor interface
serves as the absorbing boundary while other boundaries serve as reflecting boundaries due to the
tailored properties of the donor and the acceptor. As derived in \cite{Chen:2016}, such a
problem can be modeled by a diffusion-type equation with appropriate boundary conditions, which will
be introduced in \secref{sec:model}. Since the donor-acceptor interface is not exposed to the air/vacuum
and the resolution of the surface morphology is limited by the resolution of atomic force microscopy, this interface
is subject to an uncertainty with amplitude around $1\;$nm. At a first glance, this uncertainty
does not seem to affect the observation very much since its amplitude is much smaller than the film thickness.
However, in some scenarios \cite{Linetal:2013}, the fitted exciton diffusion lengths are sensitive to the uncertainty,
which may affect a chemist to determine which material should be used for a specific device.
Therefore, it is desirable to understand the quantitative effect of such an uncertainty on the
exciton diffusion length and provide a reliable estimation method  to
select appropriate models for organic materials with different crystalline orders.

Uncertainty quantification is an emerging research field that addressing these issues \cite{Xiu:09,LeMaitre:2010,Smith:2013}.
Due to the complex nature of the problems considered here, finding analytical solutions is almost impossible, so numerical methods are  very  important to study these solutions. Here we give a briefly introduction of existing numerical methods, which can be classified into non-intrusive sampling methods and intrusive methods.

Monte Carlo (MC) method is the most popular non-intrusive method \cite{glasserman:03}.
For the randomness in the partial differential equations (PDEs), one first generates $N$ random samples, and then solves the corresponding deterministic problem to obtain solution samples. Finally, one estimates the statistical information by ensemble averaging. The MC method is easy to implement, but the convergence rate is merely $O(\frac{1}{\sqrt{N}})$. Later on, quasi-Monte Carlo methods \cite{Caflisch:98} and multilevel Monte Carlo methods \cite{Giles:08} have been developed to speed up the MC method. Stochastic collocation (SC) methods explore the smoothness of PDE solutions with respect to random variables and use certain quadrature points and weights to compute solution realizations \cite{Xiu:05,Babuska:07,Webster:08}. Exponential convergence can be achieved for smooth solutions, but the quadrature points increase exponentially fast as the number of random variables increases, known as the {\it curse of dimensionality}. Sparse grids were introduced to reduce the quadrature points to some extent \cite{Griebel:04}. For high-dimensional PDEs with randomness, however, the sparse grid method is still very expensive.

In intrusive methods, solutions of the random PDEs are represented by certain basis functions, e.g., orthogonal polynomials. Typical examples are the Wiener chaos expansion (WCE) and polynomial chaos expansion (PCE) method. Then, Galerkin method is used to derive a coupled deterministic PDE system to compute the expansion coefficients. The WCE was  introduced by Wiener in \cite{Wierner:38}. However, it did not receive much attention until Cameron provided the convergence analysis in \cite{cameron:47}. In the past two decades, many efficient methods have been developed based on WCE or PCE; see \cite{Ghanem:91,Xiu:03,Xiu:2006,babuska:04,WuanHou:06} and references therein.

When dealing with relatively small input variability and outputs that do not express high nonlinearity, perturbation type methods are
most frequently used, where the random solutions are expanded via Taylor series around their mean and truncated at a certain order \cite{Matthies:2005,Dambrine2016SINUM}. Typically, at most second-order expansion is used because the resulting system of equations are typically complicated beyond the second order. An intrinsic limitation of the perturbation methods is that the magnitude of the uncertainties should be small. Similarly, one also chooses the operator expansion method to solve random PDEs. In the Neumann expansion method, we expand inverse of the stochastic operator in a Neumann series and truncate it at a certain order. This type of method often strongly depends on the underlying operator and is typically limited to static problems \cite{yamazaki:1988,Xiu:09}.

In this paper, we employ a diffusion-type equation with appropriate boundary conditions as the forward
model and the exciton diffusion length is extracted in an inverse manner. Surface roughness is
treated as a random function. After nondimensionalization, the forward model becomes a diffusion-type
equation on the domain whose geometric boundary is subject to small perturbations. Therefore, we propose
an asymptotic-based method as the forward solver with its accuracy justified both analytically and numerically.
It only requires solving several deterministic problems over the regular domain without randomness.
The efficiency of our approach is demonstrated by comparing with the SC method
as the forward solver. Of experimental interest, we find that the correlation length of
randomness is the key parameter to determine whether a 1D surrogate is sufficient
for the forward modeling. Precisely, the larger the correlation length, the more accurate the 1D surrogate.
This explains why the 1D surrogate works well for organic semiconductors with high crystalline order.

The rest of the paper is organized as follows. In \secref{sec:model}, a diffusion-type equation
is introduced as the forward model and the exciton diffusion length is extracted by solving an
inverse problem. Domain mapping method and the asymptotic-based
method are introduced in \secref{sec:method} with simulation results presented in \secref{sec:result}.
Conclusion is drawn in \secref{sec:conclusion}.

\section{Model}\label{sec:model}
\noindent
In this section, we introduce a diffusion-type equation over the random domain as the
forward model and the extraction of exciton diffusion length is done by solving an
inverse problem.

\subsection{Forward model: A diffusion-type equation over the random domain}
\begin{figure}[htbp]
\includegraphics[width=0.5\textwidth]{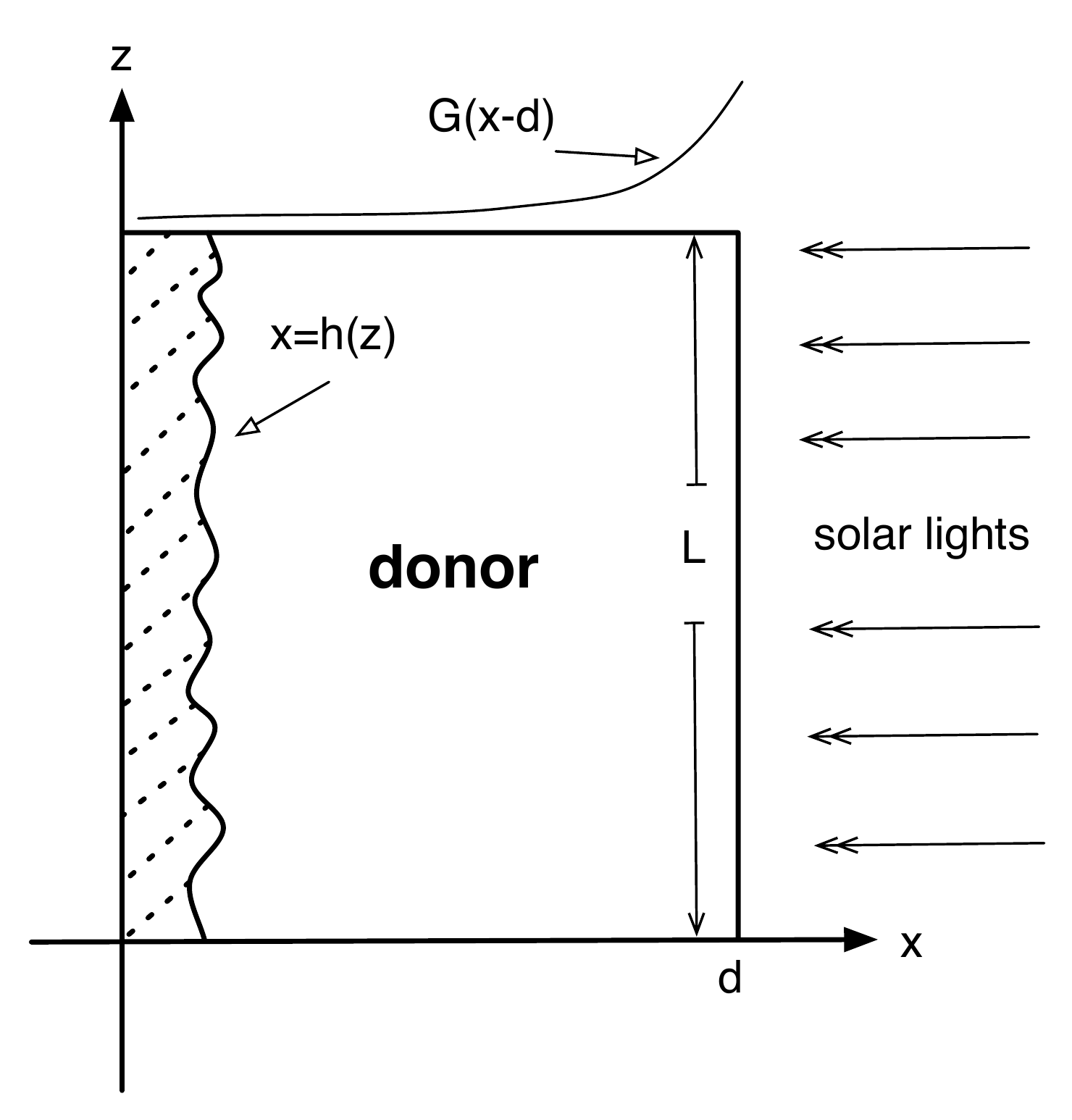}
\caption{The donor-acceptor bilayer device with film thickness $d$ along the $x$ direction and in-plane
dimension $L$ along the $z$ direction under the illumination of sun lights. One realization of the donor-acceptor
interface with uncertainty is described by $x=h(z)$. $G(x)$ is the normalized exciton generation density
which depends on $x$ only and is a decreasing function due to the phonon absorption in the donor layer.}
\label{fig:dom}
\end{figure}

\noindent
Consider a thin layer of donor located over the two dimensional domain
 $\set{(x,z): x\in (h(z,\omega),d), z\in (0,L)}$,
where $L\gg d$. Refer to Figure \ref{fig:dom}. The donnor-acceptor interface,
$\Gamma$, is described by $x=h(z,\omega)$,
a random field with period $L$:
\begin{equation}
\label{eqn:randominterface}
h(z,\omega) = \bar{h}\sum_{k=1}^{K}\lambda_{k}\theta_k(\omega)\phi_k(z),
\end{equation}
where $\{\theta_k\}$ are i.i.d. random variables, $\phi_k(z) = \sin(2k\pi\frac{z}{L})$, and $\lambda_{k}>0$ are eigenvalues that
control the decay speed of physical mode  $\phi_k(z)$.
In principle, one could  also add the cosine modes in the basis functions $\set{\phi_k}$.
We here only use the sine modes for simplicity.
In the experiment,   $\bar{h}\sim 1\;$nm due to the surface roughness limited by the
resolution of atomic force microscopy.
The thickness $d$ varys between $10\sim100\;$nm
in a series of devices.
Therefore, the   dimensionless parameter characterizing the ratio between
measurement uncertainty and film thickness
$$\eps = \bar{h}/d,$$
ranges around  $  [0.01, 0.1].$
So, it is assume that the amplitude $\bh\ll d$ in our models.
The in-plane dimensions of the donor layer are of centimeters in the experiment, but we choose
$L\sim 100\;$nm and set up the periodic boundary condition along the $z$ direction based on the
following two reasons. First, the current work treats exciton diffusion length as a homogeneous
macroscopic quantity, which is a good approximation for ordered structures.
For example, small molecules are the simplest and can form crystal structures under careful fabrication
conditions \cite{DirksenRing:1991, Rodetal:2013}. Second, the light intensity and hence the exciton
generation density is a single variable function depending on $x$ only.

Define the domain $\dom_\eps:=\set{(x,z): x\in (h(z,\omega),d), z\in (0,L)}$.
The diffusion-type equation reads as
\begin{numcases}
{ \label{eqn:rPDE2d} }
 \ldsq \left( u_{xx}(x,z)+ u_{zz}(x,z)\right) - u(x,z) + G(d-x)   = 0,
&   $(x,z)\in \dom_\eps$ \label{eqn:rPDE2d1}
\\
u_x(d, z) = 0, ~\ ~ u(h(z,\omega), z)   = 0,  &   $ 0<z<L$ \label{eqn:rPDE2d2}
\\
u(x, z)   =  u(x, z+L), &   $ h(z,\omega) < x <d$. \label{eqn:rPDE2d3}
\end{numcases}
Here $\ld$ is the exciton diffusion length which is an unknown parameter, and
the $\ldsq$ term in \eqref{eqn:rPDE2d1} describes the exciton diffusion.
Exciton-environment interaction makes some excitons emit phonons and die out,
which is described by the   term $-u$ in \eqref{eqn:rPDE2d1}. The normalized
exciton generation function $G$ is $\Real^+$-valued, and is smooth on $\Real^+\cup\set{0}$.
By solving the Maxwell equation over the layered device, one can find that $G(x)$ is a
combination of exponential functions which decay away from $0$ \cite{Born:1965}.
$x=d$ is served as the reflexive boundary and homogeneous Neumann boundary condition
is thus used there, while $x=h(z,\omega)$ is served as the absorbing boundary and homogeneous
Dirichlet boundary condition is used in \eqref{eqn:rPDE2d2}.
Periodic boundary condition is imposed along the $z$ direction in \eqref{eqn:rPDE2d3}.
It is not difficult to see that the solution $u$ to \eqref{eqn:rPDE2d}
is strictly positive in $\dom_\eps$ by the maximum principle.

The (normalized) photoluminescence is computed by the formula
\begin{equation}
\label{eqn:PL2d}
\PL[\ld,d]=\frac{1}{L}\int_{0}^L\int_{h(z,\omega)}^d u(x,z)\text{d}z\text{d}x.
\end{equation}

If the interface $\Gamma$ is random but entirely flat, i.e., $h(x,\omega)=\xi(\omega)$
for some random variable $\xi$, then the domain is a rectangle $ (\xi(\omega),d)\times(0,L)$.
Notice  that in \eqref{eqn:rPDE2d}, $G$ is a function of $x$ only.
Then, \eqref{eqn:rPDE2d} actually reduces to the following 1D problem
\begin{numcases}
{  \label{eqn:rPDE1d}}
\ldsq u_{xx}(x) - u(x) + G(d-x)  = 0,\ \ x\in (\xi,d)
\\
u_x(d) = 0, \quad u(\xi)  =  0.
\end{numcases}

For the 1D model \eqref{eqn:rPDE1d}, when $L\rightarrow 0$, the photoluminescence defined
by \eqref{eqn:PL2d} reduces to
\begin{equation}
\label{eqn:rPL1d}
\PL(\ld,d)=  \int_{\xi}^d  u(x) \ud x.
\end{equation}
This is why the normalized factor $1/L$ is used in \eqref{eqn:PL2d}.
Due to the simple analytical formula,
the 1D model given by  \eqref{eqn:rPDE1d} and \eqref{eqn:rPL1d} has been widely used to fit
experimental data for photoluminescence measurement \cite{Linetal:2013}
and photocurrent measurement \cite{Guideetal:2013}.

Since the roughness of the interface is taken into account, problem \eqref{eqn:rPDE2d}
with the random interface $\Gamma$ is viewed as a generalized and more realistic model.
The 1D model \eqref{eqn:rPDE1d} still has the uncertainty of the boundary but
fails to include the spatial variety of the donor-interface interfacial layer.
We are interested in identifying under which condition the 1D model can be viewed as
a good surrogate for the 2D model and how this condition can be related to the property
of organic semiconductors.

\subsection{Inverse problem: Extraction of exciton diffusion length}
\noindent
In the experiment, photoluminescence data $\{\wt{\PL}_i\}_{i=1}^N$ are measured for a series of
bilayer devices with different thicknesses $\{d_i\}_{i=1}^N$. Here $i$ denotes the $i$-th observation
in the experiment with $d_i$ the thickness of the donor layer. $\ld$ is the unknown parameter, and the
optimal $\ld$ is expected to reproduce the experimental data $\{d_i, \wt{\PL}_i\}_{i=1}^N$ in a proper
sense.

To achieve this, we propose the following minimization problem in the sense of mean square error
\begin{equation}
\label{min1}
\min_{\sigma} ~J(\ld) =\frac{1}{N}\sum_{i=1}^N\
\left( \mathbb{E}_\omega
[
\PL(\ld, d_i) ]-
 \wt{\PL}_i
\right)^2.
\end{equation}

We use the Newton's method to solve \eqref{min1} for $\ld$.
Given $\ld^{(0)}$, for $n=1,2,\ldots,$ until convergence, we have
\begin{equation}
\label{estimate_sigma}
\ld^{(n)} = \ld^{(n-1)} - \alpha_n \dfrac{\frac{\partial}{\partial \ld} J(\ld^{(n-1)})}
{\frac{\partial^2}{\partial \ldsq}J(\ld^{(n-1)})}.
\end{equation}
Here $\alpha_n\in (0,1]$ is given by a line search \cite{Nocedal:1999}. Details are given in Appendix \ref{sec:Newton}.

\section{Methods for solving the forward model}\label{sec:method}
\noindent
In the photoluminescence experiment, the surface roughness is very small compared to the film thickness,
i.e., $\bar{h}\sim 1\;$nm and $10\le d\le 100\;$nm. Based on this observation, we propose an
asymptotic-based method for solving the diffusion-type equation over the random domain. For comparison,
we first describe the domain mapping approach \cite{Xiu:2006}.

\subsection{Domain mapping method}\label{sec:gpc}
\noindent
To handle the random domain $\dom_\eps$, we introduce the following
transformation
\[ \tilde{y} = \frac{x - h(z,\omega)}{d-h(z,\omega)}, \quad \tilde{z} = z/L, \]
so that $\dom_\eps$ becomes  the unit square $\mathcal{D}_{\textrm{s}}=(0,1)\times(0,1)$.
Under this change of variables, Eq. \eqref{eqn:rPDE2d} becomes
the following PDE with random coefficients
(still use $y$ and $z$ to represent  $\wt{y}$ and $\wt{z}$, respectively)

\begin{equation}\label{eqn:PDE2d}
\ldsq\mathcal{L} u - u + g(y,z,\omega) =0, \quad (y, z) \in \mathcal{D}_{\textrm{s}},
\end{equation}
where the spatial differentiation operator is defined for a random element $\omega$ in the probability space
\begin{equation}\label{L}
\begin{split}
\mathcal{L}:=&
 \frac{(1-y)^2(h')^2+1}{(d-h)^2} \partial_{yy}
 +\frac{1}{L^2} \partial_{zz}
 -\frac{2}{L}\frac{(1-y)h'}{(d-h)}\partial_{yz}
\\
& -2\frac{(1-y)(h')^2}{(d-h)^2}\partial_y
-\frac{(1-y)h''}{(d-h)}\partial_y.
\end{split}
\end{equation}
and
\begin{equation}\label{g}
g(y,z,\omega):=G((1-y)(d-h(z,\omega))).
\end{equation}

The boundary condition is
\begin{equation}\label{eqn:bc2d}
\begin{split}
 \partial_y u(1,z)  =   0, \quad u(0,z)= 0, \quad z\in (0,1), \\
u(y,z) =   u(y,z+1),  \quad y\in (0,1).
\end{split}
\end{equation}
The photoluminescence defined in \eqref{eqn:PL2d} is then transformed into
\begin{equation}
\label{eqn:PL2D}
\PL(\ld, d)= \int_{0}^1\int_{0}^1 u(y,z)(d-h(z,\omega))\text{d}y\text{d}z.
\end{equation}

\begin{remark}
In 1D, changing of variable $y = \frac{x - \xi}{d-\xi}$ also transforms
\eqref{eqn:rPDE1d} to
a differential equation with random coefficients over the unit interval.
\begin{equation}\label{eqn:PDE1d}
\ldsq\mathcal{L}_1 u(y) - u(y) + G((1-y)(d-\xi)) =  0,  \quad y\in (0,1)
\end{equation}
with
\begin{equation}\label{L1d}
\mathcal{L}_1:= \frac{1}{(d-\xi)^2} d_{yy}
\end{equation}
and the boundary condition
\begin{equation}\label{eqn:bc1d}
u_y(1) = 0, \quad u(0) =  0.
\end{equation}

Accordingly, the photoluminescence can be written as
\begin{equation}
\label{eqn:PL1d}
\PL(\ld, d)=(d-\xi)\int_{0}^1 u(y)\text{d}y.
\end{equation}
\end{remark}

\begin{remark}
The generation term in \eqref{g} depends on both $y$ and $z$
after changing of variables. We expect some dimensional effect on the
estimation of $\ld$, which will be carefully examined in \secref{sec:result}.
\end{remark}

\subsection{Finite difference method for the model problem}
\noindent
We use finite difference method to discretize the forward model \eqref{eqn:PDE2d} developed in \secref{sec:gpc}. We partition the domain $\mathcal{D}_{\textrm{s}}=[0,1]\times[0,1]$ into $(N_y+1)\times (N_z+1)$ grids with meshes $h_y=\frac{1}{N_y}$ and $h_z=\frac{1}{N_z}$.
Denote by $u_{i,j}$ the numerical approximation of $u(y_i,z_j)$, where $y_i=(i-1)h_y$, $z_j=(j-1)h_z$ with $i=1,...,N_y+1$ and $j=1,...,N_z+1$, respectively. For the discretization in space, we use a second-order, centered-difference scheme \cite{morton:2005}. We introduce the difference operators
\[
D_{0}^{y}u_{i,j}=\frac{u_{i+1,j}-u_{i-1,j}}{2h_y},  \quad D_{-}^{y}u_{i,j}=\frac{u_{i,j}-u_{i-1,j}}{h_y}, \quad
D_{+}^{y}u_{i,j}=\frac{u_{i+1,j}-u_{i,j}}{h_y}.
\]
The operators $D_{0}^{z}$, $D_{-}^{z}$, and $D_{+}^{z}$ are defined similarly.
For each $\omega \in \Omega$ and each interior mesh point $(i,j)$ with $2 \leq i \leq N_y,  2 \leq j \leq N_z$,
we discretize the forward model \eqref{eqn:PDE2d} as
\begin{align}
&\sigma^2 \frac{(1-y_i)^2(h')^2+1}{(d-h)^2}D_+^y D_-^y u_{i,j}
  + \frac{\sigma^2 }{L^2}D_+^z D_-^z u_{i,j} - \frac{2\sigma^2 }{L}\frac{(1-y_i)h'}{(d-h)}D_0^y D_0^z u_{i,j}\nn \\
& - \left(2\sigma^2 \frac{(1-y_i)(h')^2}{(d-h)^2} + \sigma^2 \frac{(1-y_i)h''}{(d-h)} \right)  D_0^y u_{i,j} - u_{i,j}
 =- g(y_i,z_j,\omega),
\label{eqn:PDE2d_FDM}
\end{align}
where $h$, $h'$, and $h''$ are evaluated at $(y_i,z_j)$.

We then discretize the boundary conditions \eqref{eqn:bc2d} on $\partial\mathcal{D}_{\textrm{s}}$. The Dirichlet boundary condition
on $y=0$ gives $u_{1,j}=0$,  $1 \leq j \leq N_z+1$. For the Neumann boundary condition
on $y=1$, we introduce ghost nodes at $(y_{-1},z_{j})$ and obtain a second order accurate finite difference approximation
$\frac{u_{1,j}-u_{-1,j}}{2h_y}=0$. Then, the values of the $u_{-1,j}$ at the ghosts nodes are eliminated by combining with
Eq. \eqref{eqn:PDE2d_FDM}. Finally, the periodic boundary condition along the $z$ direction gives $u_{i,N_z+1}=u_{i,1}$.
We solve a system of $N_y(N_z+1)$ linear equations for $\{u_{i,j}\}$ with
$2 \leq i \leq N_y+1$ and $1 \leq j \leq N_z+1$.

The equations have a regular structure, each equation involving at most nine unknowns. Thus the corresponding matrix of the system is sparse and can be solved efficiently using existing numerical solvers. After obtaining $\{u_{i,j}\}$, we use the 2D trapezoidal quadrature rule to compute the photoluminescence $\PL(\ld, d)$ defined in \eqref{eqn:PL2D}.

In this paper, we choose the sparse-grid based SC method \cite{Griebel:04,Webster:08} to discretize the stochastic dimension
in Eq. \eqref{eqn:PDE2d}. As such the expectation of $u(y,z,\omega)$ is computed by
\[
\mathbb{E}[u(y,z,\omega)] = \sum_{q=1}^{Q} u(y,z,s_q)w_q, \label{eqn:sparse_grid}
\]
where $s_q$ are sparse-grid quadrature points, $w_q$ are the corresponding weights, and $Q$ is the number of sparse-grid points.
Other functionals of $u(y,z,\omega)$ can be computed in the same way.  When the solution $u(y,z,\omega)$ is smooth in the stochastic dimension, the SC method provides very accurate results.

\subsection{An asymptotic-based method}\label{sec:asymptotic}
\noindent
If we rewrite Eq. \eqref{eqn:rPDE2d} in the nondimensionalized form with the change of variables $\tilde{x} = x/d$ and $\tilde{z} = z/L$,
the domain $\dom_\eps$ becomes
\[
\dom_{\textrm{s}, \eps}:=\set{(x,z) \in (\eps\tilde{h}(z,w),1)\times(0,1)},
\]
where $\eps= \bar{h}/d$. When $\eps = 0$, $\dom_{\textrm{s}, \eps}$ becomes $\dom_{\textrm{s}, 0}=\dom_{\textrm{s}}=(0,1)\times(0,1)$.
Here
\begin{equation}\label{eqn:htilde}
\tilde{h}(z,w) = \sum_{k=1}^{K} \lambda_k\theta_k(\omega)\phi_k(z),
\end{equation}
where $K$ is the mode number in the interface modeling. As discussed in \secref{sec:model}, $\eps\sim 0.01-0.1$.
Therefore, it is meaningful to derive the asymptotic equations when $\eps\rightarrow 0$.
For ease of description, we list the main results below. The main idea is: (1) we rewrite Eq. \eqref{eqn:rPDE2d}
over $\dom_{\textrm{s}, \eps}$; (2) with appropriate extension/restriction
of solutions on the fixed domain $\dom_{\textrm{s}}$, we obtain a Taylor series with each term satisfying a PDE of
the same type with the boundary condition involving lower order terms; (3) we apply the inverse transform for each
term and change the domain $\dom_{\textrm{s}}$ back to $\dom_0=(0,d)\times(0,L)$.
Detailed derivation can be found in Appendix \ref{sec:derivation} for self-consistency.
The interested readers can find the systematic study on asymptotic expansions for more general problems in \cite{Chen:2017}.

The asymptotic expansion over the fixed domain $\dom_0$ is of the form
\begin{equation}
\label{series}
w_\eps(x,z)=\sum_{n=0}^\infty \eps^n w_{n}(x,z) \quad \text{for }(x,z)\in \overline{\dom_0}.
\end{equation}
The equation for each $w_n$ can be derived in a sequential manner. Only the first three
terms are listed here. More details are included in Appendix \ref{sec:derivation}.

The leading term $w_0(x,z)$ is the solution to the boundary value problem
\begin{equation}\label{eqn:w0}
\begin{cases}
&\sigma^2\partial_{xx}w_0+\sigma^2\partial_{zz}w_0-w_0
+G(d-x)=0 \quad \text{in } \dom_0,\\
&\partial_x w_0(d,z)=0, \\
&w_{0}(0,z)=0,
\quad  \text{for }  0\leqslant z\leqslant L,\\
&w_0(x,z+L)=w_0(x,z), \quad \text{for }  0\leqslant x\leqslant d,
\end{cases}
\end{equation}
and $w_1(x,z,\omega)$ solves
\begin{equation}\label{eqn:w1}
\begin{cases}
&\sigma^2\partial_{xx}w_1+\sigma^2\partial_{zz}w_1-w_1
=0 \quad \text{in } \dom_0,\\
&\partial_x w_1(d,z,\omega)=0, \\
&w_{1}(0,z,\omega)=-d\tilde{h}(z,\omega)\partial_x w_{0}(0,z),
\quad  \text{for }  0\leqslant z\leqslant L,\\
&w_1(x,z+L,\omega)=w_1(x,z,\omega), \quad \text{for }  0\leqslant x\leqslant d.
\end{cases}
\end{equation}
$w_2(x,z,\omega)$ is the solution to the following boundary value problem
\begin{equation}\label{eqn:w2}
\begin{cases}
&\sigma^2\partial_{xx}w_2+\sigma^2\partial_{zz}w_2-w_2
=0 \quad \text{in } \dom_0,\\
&\partial_x w_2(d,z,\omega)=0, \quad \text{for } 0\leqslant z\leqslant L,\\
&w_{2}(0,z,\omega)=-d\tilde{h}(z,\omega)\partial_x w_{1}(0,z,\omega)+\frac{(d\tilde{h}(z,\omega))^2}{2\sigma^2}G(d),
\   \text{for }  0\leqslant z\leqslant L,\\
&w_2(x,z+L,\omega)=w_2(x,z,\omega), \quad \text{for }  0\leqslant x\leqslant d.
\end{cases}
\end{equation}
\begin{remark}
As demonstrated in Eqs. \eqref{eqn:w0}, \eqref{eqn:w1}, and \eqref{eqn:w2}, the asymptotic expansion in \eqref{series}
requires a sequential construction from lower order terms to high order terms and the partial derivatives of lower terms appear
in the boundary condition for high terms. Numerically, we use the second-order
finite difference scheme for \eqref{eqn:w0}, \eqref{eqn:w1}, and \eqref{eqn:w2}.
For boundary conditions, we use the one-sided beam warming scheme to discretize $\partial_x w_{0}(0,z)$ and
$\partial_x w_{1}(0,z,\omega)$ so the overall numerical schemes are still of second order accuracy.
\end{remark}

Define $v^{[n]}=\sum_{k=0}^n \eps^k w_k$.
Note that $w_0$ is a function of $(x,z)$ only. The zeroth order approximation of the photoluminescence is
\begin{equation}
\label{I0}
\PL[u]\approx\PL[v^{[0]}] = \frac1L\int_{\dom_\eps} w_0(x,z)\,\rd x\rd z\approx\frac1L\int_{\dom_0} w_0(x,z)\,\rd x\rd z=:\PL_0[v^{[0]}],\\
\end{equation}
and so
\begin{equation}
\label{EI0}
\EE[\PL[u]]\approx\EE\left[\PL_0[v^{[0]}]\right]=\PL_0[w_0].
\end{equation}

For $k=1,2,\ldots, K$, let $w_{1,k}(x,z)$ be the solution to \eqref{eqn:w1} with $\phi_k(z)$ in place of $\tilde{h}(z,\omega)$
\begin{equation}\label{eqn:w1k}
\begin{cases}
&\sigma^2\partial_{xx}w_{1,k}+\sigma^2\partial_{zz}w_{1,k}-w_{1,k}
=0 \quad \text{in } \dom_0,\\
&\partial_x w_{1,k}(d,z)=0, \\
&w_{1,k}(0,z)=-d\phi_k(z)\partial_x w_{0}(0,z),
\quad  \text{for }  0\leqslant z\leqslant l,\\
&w_{1,k}(x,z+L )=w_{1,k}(x,z), \quad \text{for }  0\leqslant x\leqslant d.
\end{cases}
\end{equation}
Then by linearity, the solution $w_1$ to \eqref{eqn:w1} with $\tilde{h}$ given by \eqref{eqn:htilde}
can be expressed as
\begin{equation}\label{eqn:KLu1}
w_1(x,z,\omega)=\sum_{k=1}^K \lambda_k\theta_k(\omega)w_{1,k}(x,z).
\end{equation}
 Hence the first order approximation of the photoluminescence becomes
\begin{equation}
\label{I1}
\begin{split}
\PL[u]&\approx\PL[v^{[1]}]=\frac1L\int_{\dom_\eps}v^{[1]}(x,z,\omega)\,\rd x\rd z\approx\frac1L\int_{\dom_0}v^{[1]}(x,z,\omega)\,\rd x\rd z\\
&=\frac1L\int_{\dom_0} \left[w_0(x,z)
+\eps w_1(x,z,\omega)\right]\,\rd x\rd z \\
&=\frac1L\int_{\dom_0} w_{0}(x,z)\,\rd x\rd z+
\frac \eps L\sum_{k=1}^K \lambda_k\theta_k(\omega)\int_{\dom_0} w_{1,k}(x,z)\,\rd x\rd z \\
&=\PL_0[w_0]+\eps\sum_{k=1}^K \lambda_k\theta_k(\omega) \PL_0[w_{1,k}]
=:\PL_1[v^{[1]}],
\end{split}
\end{equation}
and so
\begin{equation}
\label{EI1}
\EE[\PL[u]]\approx\EE\left[\PL_1[v^{[1]}]\right]
=\PL_0[w_0]+\eps\sum_{k=1}^K\lambda_k\EE[\theta_k] \PL_0[w_{1,k}].
\end{equation}

Next, we consider the second order approximation of the photoluminescence.
Since $\tilde{h}$ and $w_1$ are given by \eqref{eqn:htilde} and \eqref{eqn:KLu1},
the boundary condition for $w_2$ at $x=0$ can be written as
\[
w_2=\sum_{j,k=1}^K \lambda_j\lambda_k\theta_j\theta_k\left(-d\phi_j\partial_xw_{1,k}+\frac{G(d)}{2\sigma^2}d^2\phi_j\phi_k\right).
\]
Introduce $w_{2,j,k}(x,z)$ as the solution to the boundary value problem \eqref{eqn:w2} with
the boundary condition at $x=0$ replaced by
\[
w_{2,j,k}=-d\phi_j\partial_xw_{1,k}+\frac{G(d)}{2\sigma^2}d^2\phi_j\phi_k.
\]
Then
\[
w_2(x,z,\omega)=\sum_{j,k=1}^K \lambda_j\lambda_k\theta_j(\omega)\theta_k(\omega) w_{2,j,k}(x,z),
\]
and consequently, the second order approximation of the photoluminescence is
\begin{equation}
\label{I2}
\begin{split}
\PL[u]\approx~&\PL[v^{[2]}]=\frac1L\int_{\dom_\eps} v^{[2]}(x,z,\omega)\,\rd x\rd z\\
\approx~&\frac1L\int_{\dom_0} v^{[2]}(x,z,\omega)\,\rd x\rd z-\frac \eps{2L}\int_0^L v^{[2]}(0,z,\omega)h(z,\omega)\,\rd z\\
\approx~&\frac1L\int_{\dom_0} [w_0+\eps w_1+\eps^2 w_2]\,\rd x\rd z
-\frac \eps{2L}\int_0^L [w_0+\eps w_1](0,z,\omega)h(z,\omega)\,\rd z\\
=~&\frac1L\int_{\dom_0} [w_0+\eps w_1+\eps^2 w_2]\,\rd x\rd z+\frac {\eps^2}{2L}\int_0^L
[d\tilde{h}(z,\omega)]^2\partial_x w_0(0,z)\,\rd z\\
=~&\PL_0[w_0]+\eps\sum_{k=1}^K\lambda_k\theta_k(\omega) \PL_0[w_{1,k}]+\eps^2\sum_{j,k=1}^K
\lambda_j\lambda_k\theta_j(\omega) \theta_k(\omega)  \PL_0[w_{2,j,k}]\\
&+\frac {\eps^2d^2}{2L}\sum_{j,k=1}^K \lambda_j\lambda_k\theta_j(\omega)\theta_k(\omega)\int_0^L \phi_j(z)\phi_k(z)\partial_x w_0(0,z)\,\rd z\\
=:&\PL_2[v^{[2]}],
\end{split}
\end{equation}
and we have
\begin{equation}
\label{EI2}
\begin{split}
\EE[\PL[u]]\approx\EE\left[\PL_2[v^{[2]}]\right]=&\PL_0[w_0]+
\eps\sum_{k=1}^K\lambda_k\EE[\theta_k] \PL_0[w_{1,k}]+\eps^2\sum_{j,k=1}^K
\lambda_j\lambda_k\EE[\theta_j\theta_k] \PL_0[w_{2,j,k}]\\
&+\frac {\eps^2d^2}{2L}\sum_{j,k=1}^K \lambda_j\lambda_k\EE[\theta_j\theta_k]\int_0^L \phi_j(z)\phi_k(z)\partial_x w_0(0,z)\,\rd z.
\end{split}
\end{equation}

In general, $w_n$ can be written as the sum of $K^n$ functions, each of which solves a deterministic problem.

The approximation accuracy of a finite series in \eqref{series} is given by the following theorem.
Proof can be found in \cite{Chen:2017}.
\begin{theorem}\label{theorem}
Assume $\dom_0\subset\dom_\eps\subset\dom_{\eps_0}$ with $\eps\in [0,\eps_0]$ and $\partial\dom_0\in C^{\infty}$.
Also assume $G\in C^{\infty}(\overline{\dom_{\eps_0}})$ and $h\in C^{\infty}(\partial\dom_0)$.
Then, $\forall n, m\geq 0$,
\begin{equation}\label{eqn:vn err}
\bigl\| v^{[n]}(\omega) -u(\omega) \bigr\|_{H^m(\dom_0)} = \Od (\eps^{n+1}) \quad \mathbb{P}-\textrm{a.e.\ } \omega\in\Omega,
\end{equation}
where $u$ is the solution to \eqref{eqn:rPDE2d} and $v^{[n]}=\sum_{k=0}^n \eps^k w_k$.
\end{theorem}
To proceed, let us recall the definition of Bochner spaces.
\begin{definition}
Given a real number $p\geq 1$ and a Banach space $X$, the Bochner space is
\[
L_{\mathbb{P}}^p(\Omega,X) = \{u:\Omega\rightarrow X \; | \; \norm{u}_{L_{\mathbb{P}}^p(\Omega,X)} \textrm{is finite}\}
\]
with
\[
\norm{u}_{L_{\mathbb{P}}^p(\Omega,X)} : =
\begin{cases}
&\left( \int_{\Omega} \norm{u(\cdot,\omega)}_X^p \mathrm{d}\;\mathbb{P}(\omega)\right)^{1/p}, \quad p<\infty \\
&\mathrm{ess\;sup}_{\omega\in\Omega} \norm{u(\cdot,\omega)}_X, \quad p=\infty.
\end{cases}
\]
\end{definition}
\begin{proposition}\label{proposition}
Given $h\in L_{\mathbb{P}}^\infty(\Omega,C^1(\partial \dom_0))$, then $w_n,\; n\geq 0$ belongs to $L_{\mathbb{P}}^2(\Omega,H^1(\dom_0))$
and hence
\[
\bigl\| v^{[n]} -u \bigr\|_{L_{\mathbb{P}}^2(\Omega,H^1(\dom_0))} = \Od (\eps^{n+1}).
\]
\end{proposition}
\begin{proof}
From Theorem \ref{theorem}, for $m=1$, we have
\[
\bigl\| v^{[n]}(\omega) -u(\omega) \bigr\|_{H^1(\dom_0)} = \Od (\eps^{n+1}) \quad \mathbb{P}-\textrm{a.e.\ } \omega\in\Omega.
\]
Since $w_n,\; n\geq 0$ satisfies the same elliptic equation \eqref{eqn:w0} with a boundary condition depending on $w_k, \;k\leq n-1$.
By the Lax-Milgram's theorem, we have $w_n \in L_{\mathbb{P}}^2(\Omega,H^1(\dom_0)),\; n\geq 0$.
Therefore, $v^{[n]}\in L_{\mathbb{P}}^2(\Omega,H^1(\dom_0)),\; n\geq 0$ and the desired result is obtained.
\end{proof}
A direct consequence of Proposition \ref{proposition} is
\begin{equation}
\label{expection}
\norm{\mathbb{E}(v^{[n]}) - \mathbb{E}(u)}_{H^1(\dom_0)} = \Od (\eps^{n+1}).
\end{equation}
Based on the above assertions, we have
\begin{corollary}
\label{corollary}
For \eqref{EI0}, \eqref{EI1}, and \eqref{EI2}, we have the following approximation errors
\begin{eqnarray}
\label{expection0}
\left\lvert\EE\left[\PL_0[v^{[0]}]\right] - \EE\left[\PL[u]\right]\right\rvert = \Od (\eps^{1}), \\
\label{expection1}
\left\lvert\EE\left[\PL_1[v^{[1]}]\right] - \EE\left[\PL[u]\right]\right\rvert = \Od (\eps^{2}), \\
\label{expection2}
\left\lvert\EE\left[\PL_2[v^{[2]}]\right] - \EE\left[\PL[u]\right]\right\rvert = \Od (\eps^{3}).
\end{eqnarray}
\end{corollary}
In summary, by using the asymptotic expansion solution, we circumvent the difficulty of sampling the random function
and solving PDEs on irregular domains for each sample. In our approach, there is no statistical error or errors from
numerical quadratures as in MC method, SC method, and PCE method. However, our method is applicable only for small
perturbation of the random interface, where a small $n$ is sufficient in practice. The computational cost depends on the
approximation order $n$ and the number of modes $K$ used to represent the random interface, and increases proportionally
to $K^n$.

\section{Numerical Results}\label{sec:result}
\noindent
In this section, we numerically investigate the accuracy and efficiency of the asymptotic-based method in computing photoluminescence
and the efficiency in estimating the exciton diffusion length. In addition, we study of the validation of the diffusion-type model,
i.e., under which condition the 1D model can be viewed as a good surrogate for the 2D model.

\subsection{Accuracy and efficiency of the asymptotic-based method}
Consider the forward model defined by Eq. \eqref{eqn:rPDE2d} over
$\dom_\eps:=\set{(x,z): x\in (h(z,\omega),d), z\in (0,L)}$.
Recall that the random interface $h(z,\omega)$ between the donor and the acceptor is parameterized by
$h(z,\omega) = \bar{h}\sum_{k=1}^{K}  \lambda_{k} \theta_k(\omega)\sin(2k\pi\frac{z}{L})$,
where $\theta_k(\omega)$ are i.i.d. uniform random variables and $K$ is the
number of random variables in the model.

We first solve \eqref{eqn:PDE2d} over the fixed domain $\dom_{\textrm{s}}=(0,1)\times(0,1)$ in the domain mapping method
using the SC method. Note that the spatial differentiation operator in \eqref{L} depends on the random variables
in a highly nonlinear fashion, which makes the WCE method and PCE method extremely difficult. In the asymptotic-based method, we solve \textit{deterministic} boundary value problems
\eqref{eqn:w0}, \eqref{eqn:w1}, and \eqref{eqn:w2} over the \textit{fixed} domain $\mathcal{D}_0=(0,d)\times(0,L)$, respectively. Recall that in the asymptotic-based method, $\eps= \bar{h}/d$ and the random interface becomes $\tilde{h}(z,\omega)=\sum_{k=1}^{K} \lambda_k\theta_k(\omega)\sin(2k\pi z)$.
In our simulation, the random interface $h(z,\omega)$ is parameterized by $K=5$ random variables.
The accuracy of the asymptotic-based method is verified by two numerical tests.  In the first test,
$\theta_k\sim U(0,1)$, while in the second one $\theta_k\sim U(-1,1)$.

To compute the reference solution, we employ the finite difference method to discretize the spatial dimension of Eq. \eqref{eqn:PDE2d} with a mesh size $H=\frac{1}{128}$, and use the sparse-grid based SC method to discretize the stochastic dimension. We choose level six sparse grids with 903 quadrature points. After obtaining solutions at all quadrature points, we compute
the expectation of the photoluminescence, which provides a very accurate reference solution.
In the asymptotic-based method, we use the finite difference method to discretize the spatial dimension of boundary value problems \eqref{eqn:w0}, \eqref{eqn:w1k}, and \eqref{eqn:w2} for $w_{2,j,k}$ with a mesh size $H=\frac{1}{64}$. Expectations
$\EE[\theta_k]$ in \eqref{EI1} and  $\EE[\theta_j\theta_k]$ in \eqref{EI2} can be easily computed beforehand. Therefore,
given the approximate solutions $w_0$, $w_{1,k}$, and $w_{2,j,k}$, we immediately obtain different order approximations of the expectation of the photoluminescence. This provides the significant computational saving over the SC method.

\begin{figure}[h]
\centering
\subfigure[$U(0,1)$]{\label{fig:EstEDL_EX1_AccuAsympMethodFB4U01a}
\includegraphics[width=0.49\linewidth]{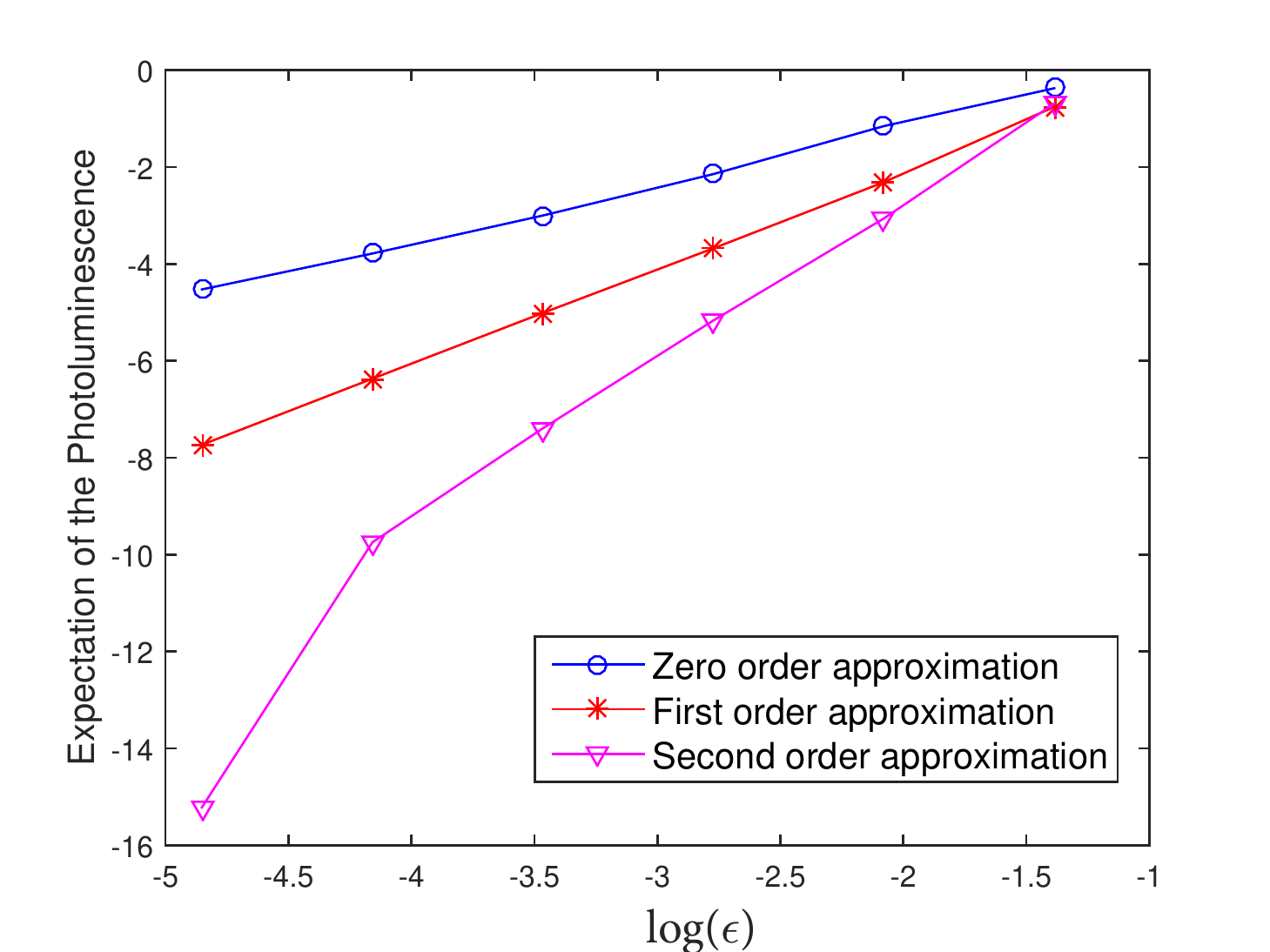}}%
\subfigure[$U(-1,1)$]{\label{fig:EstEDL_EX1_AccuAsympMethodFB4U01b}
\includegraphics[width=0.49\linewidth]{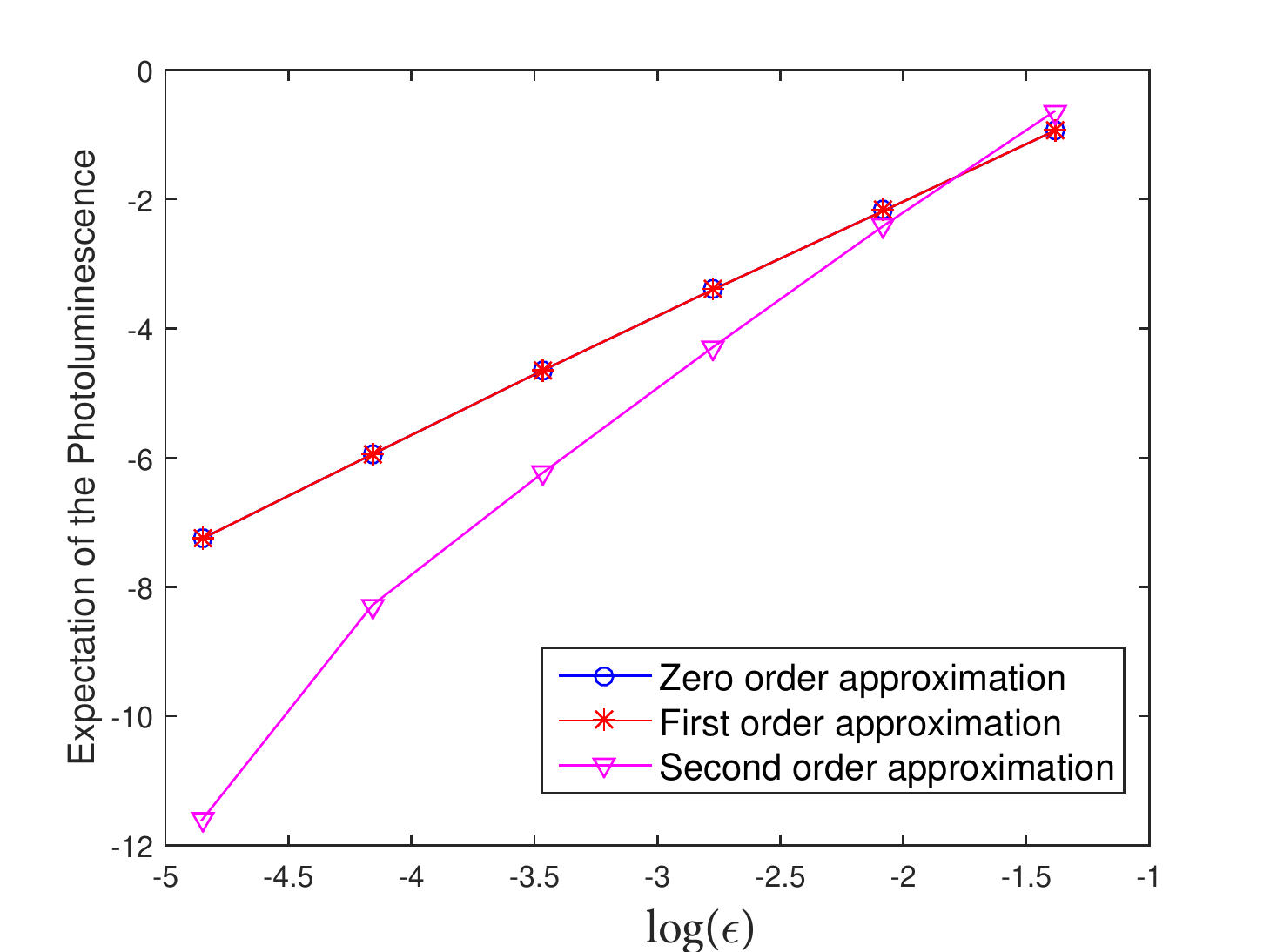}}%
\caption{\small  Convergent results of the asymptotic-based method with the zeroth, first, and second
order approximations. (a) $\theta_k \sim U(0,1)$. The slopes of the zeroth, first and
second order approximation results are 1.21, 1.99, and 3.81, respectively; (b) $\theta_k \sim U(-1,1)$. The slopes of the zeroth,
first and second order approximation results are 1.82, 1.82, and 3.06, respectively.}\label{fig:EstEDL_EX1_AccuAsympMethodFB4U01}
\end{figure}

For $\eps=2^{-i}$, $i=2,...,7$, Figure \ref{fig:EstEDL_EX1_AccuAsympMethodFB4U01} shows the approximation accuracy of the asymptotic-based method.
In Figure \ref{fig:EstEDL_EX1_AccuAsympMethodFB4U01a}, $\theta_k \sim U(0,1)$. The approximated expectation of the photoluminescence obtained by
using the zeroth, first and second order approximations are shown in the lines with circle, star, and triangle, with convergence rates
1.21, 1.99, and 3.81, respectively. In Figure \ref{fig:EstEDL_EX1_AccuAsympMethodFB4U01b}, $\theta_k\sim U(-1,1)$. In this case, $\EE[\theta_k]=0$,
so the zeroth and first order approximations produce the same results. The second order approximation provides a better result. The
corresponding convergence rates are 1.82, 1.82, and 3.06, respectively. These results confirm the theoretical estimates in Corollary \ref{corollary}.

We conclude this subsection with a discussion on the computational time of our method. In these two tests, on average it takes
164.5 seconds to compute one reference expectation of the photoluminescence. If we choose a low level SC method to compute the expectation of the photoluminescence, it takes 27.3 seconds to compute one reference expectation of the photoluminescence that gives a comparable approximation result to our asymptotic-based method. However, our method with the second order approximation only takes 1.56 seconds to obtain one result. We achieve a $18X$ speedup over the SC method. Generally, the ratio of the speedup is problem-dependent. It is expected that higher ratio of speedup can be achieved it one solves a problem where the random interface is parameterized by high-dimensional random variables.

\subsection{Estimation of the exciton diffusion length}
In this section, we estimate the exciton diffusion length in an inverse manner with the asymptotic-based method as the
forward solver. Since only limited photoluminescence data from experiments are available, we solve
the forward model \eqref{eqn:rPDE2d} to generate data in our numerical tests. Specifically, given the exciton diffusion
length $\sigma$, the exciton generation function $G$, the in-plane dimension $L$, and the parametrization of the random interface
$h(z,\omega)$, we solve Eq. \eqref{eqn:rPDE2d} for a series of thicknesses $\{d_i\}$, and calculate the corresponding
expectations of the photoluminescence data $p\{\tilde{\PL}_i\}$ according to Eq. \eqref{eqn:PL2d}.
Therefore, $\{ d_i, \tilde{\PL}_i\}$ serves as the ``experimental'' data.
We then solve the minimization problem \eqref{min1} based on our numerically generated data $\{ d_i, \tilde{\PL}_i\}$ to estimate
the ``exact'' exciton diffusion length $\sigma$ in the presence of randomness, denoted by $\sigma_{exact}$ and will be used for
comparison later.

\begin{figure}[h]
\centering
\subfigure[$\sigma=5$]{\label{fig:EstEDL_EX2_DiffusionLengthSigmaCase1}
\includegraphics[width=0.3\linewidth]{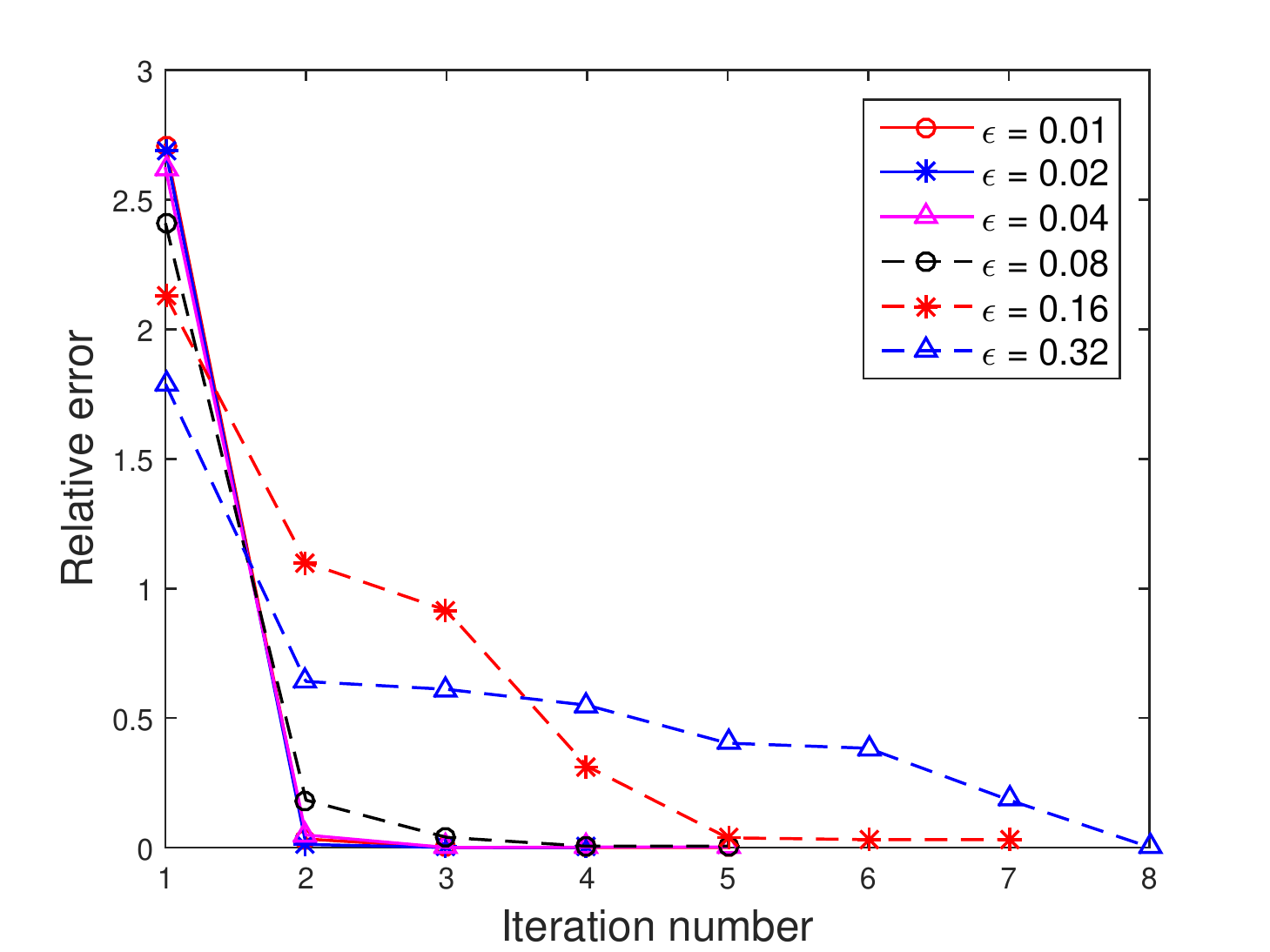}}%
\subfigure[$\sigma=10$]{\label{fig:EstEDL_EX2_DiffusionLengthSigmaCase2}
\includegraphics[width=0.3\linewidth]{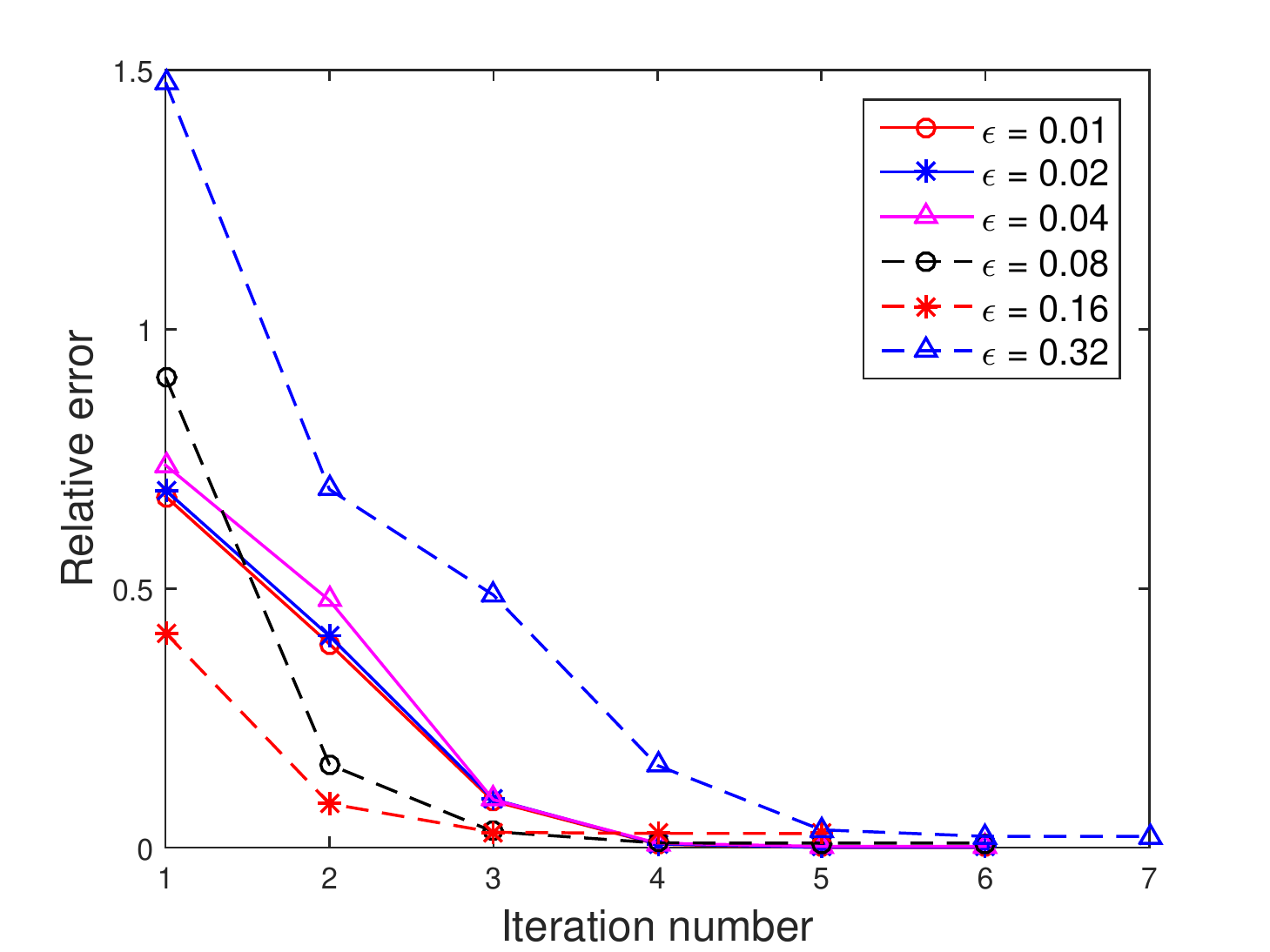}}%
\subfigure[$\sigma=20$]{\label{fig:EstEDL_EX2_DiffusionLengthSigmaCase3}
\includegraphics[width=0.3\linewidth]{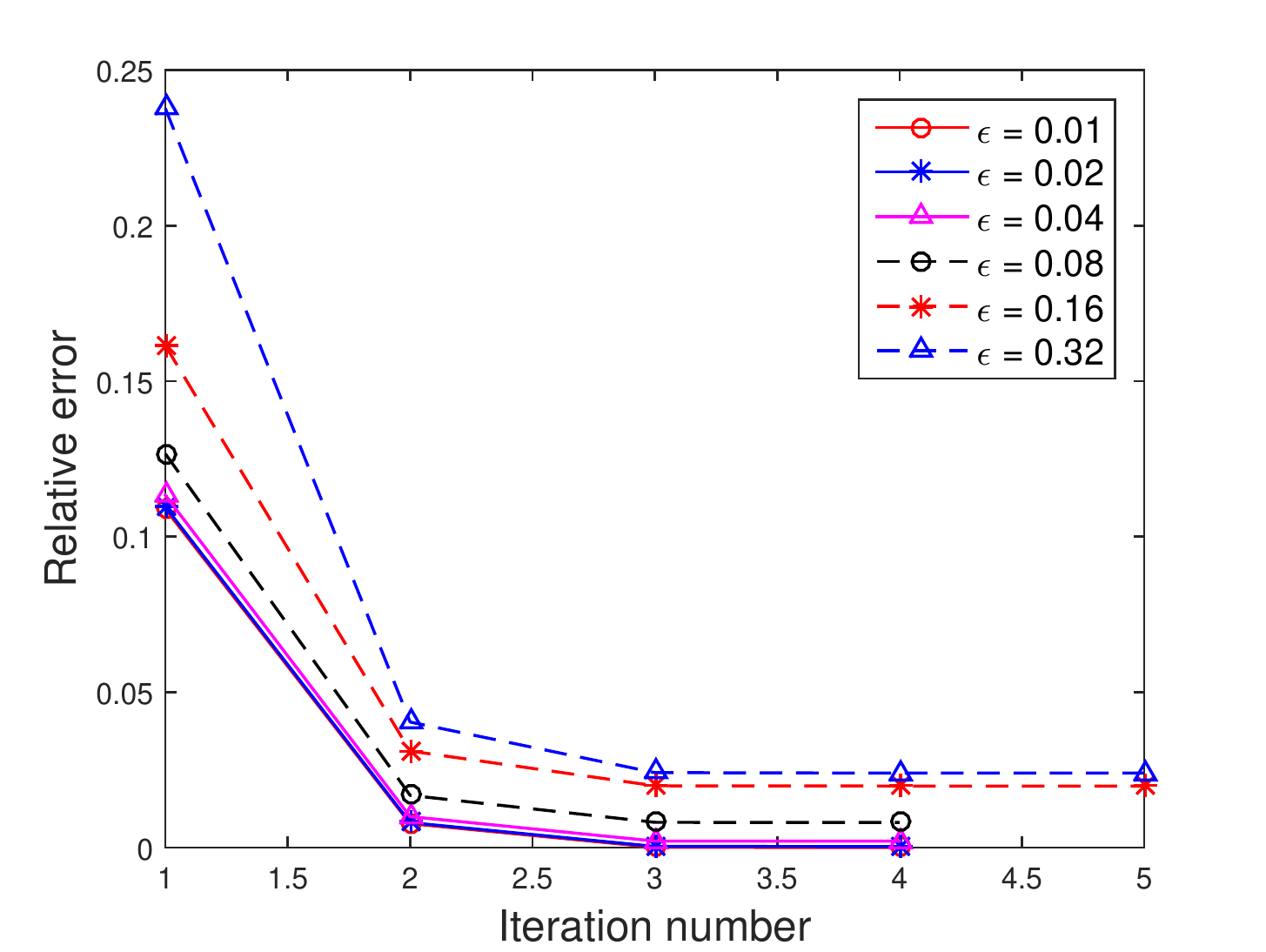}}
\caption{\small Convergence history of the exciton diffusion length for various $\eps$, measured in
the relative error defined as $E^{n,\eps}=|\frac{\sigma_{exact}- \sigma^{n,\eps}}{\sigma_{exact}}|$
with $n$ the iteration number. The ``exact'' data is obtained by the 2D model (Eqs. \eqref{eqn:rPDE2d} and \eqref{eqn:PL2d}) with
a prescribed $\sigma$. (a) $\sigma=5$; (b) $\sigma=10$; (c) $\sigma=20$.}\label{fig:EstEDL_EX2_DiffusionLengthSigmaCase}
\end{figure}

We fix $L=4$ in all our numerical tests since it is found that this minimizer is not sensitive to the in-plane dimension $L$.
We show the convergence history of exciton diffusion lengths for various $\eps$ in Figure \ref{fig:EstEDL_EX2_DiffusionLengthSigmaCase},
where the photoluminescence data are generated with $\sigma=5$, $\sigma=10$, and $\sigma=20$, respectively.
Here the relative error is defined as $E^{n,\eps}=|\frac{\sigma_{exact}- \sigma^{n,\eps}}{\sigma_{exact}}|$, where $n$ is the iteration number, $\sigma_{exact}$ is the ``exact''  exciton diffusion length, and $\sigma^{n,\eps}$ is the numerical result defined in Eq. \eqref{estimate_sigma}.
To show more details about the accuracy of our asymptotic-based method, in Tables \ref{ConvergenceOfDiffusionLengthSigmaCase1}, \ref{ConvergenceOfDiffusionLengthSigmaCase2}, and \ref{ConvergenceOfDiffusionLengthSigmaCase3}, we list the relative errors of our method for plotting Figures \ref{fig:EstEDL_EX2_DiffusionLengthSigmaCase1}, \ref{fig:EstEDL_EX2_DiffusionLengthSigmaCase2}, and \ref{fig:EstEDL_EX2_DiffusionLengthSigmaCase3}. In all numerical tests, we choose the same termination criteria $|\ld^{(n)}-\ld^{(n-1)}|<10^{-4}$ in the Newton's method.
Our asymptotic-based method performs well in estimating the exciton diffusion length. In general, the smaller amplitudes the random interface, the more accurate the exciton diffusion length and the smaller the iteration number. Additionally, for larger exciton diffusion lengths $\sigma_{exact}$, a faster convergence in the optimization approach is observed.

\begin{table}[h!]
\centering
 \begin{tabular}{|c |c |c |c |c |c |c| }
 \hline
 $n$ & $ \eps = 0.01 $ & $\eps = 0.02  $ &$\eps = 0.04  $ &$\eps = 0.08  $ & $\eps = 0.16  $ & $\eps = 0.32  $ \\
 \hline
1 &  2.711349  &  2.691114  &  2.620542  &  2.409153  &  2.132821 &  1.789203  \\
2 &  0.033009  &  0.011998  &  0.048469  &  0.182973  &  1.100014 &  0.640850  \\
3 &  0.000480  &  0.000238  &  0.000147  &  0.039389  &  0.915513 &  0.610645 \\
4 &  0.000033  &  0.000318  &  0.001678  &  0.005017  &  0.313381 &  0.549289  \\
5 &  0.000034  &  0.000317  &  0.001679  &  0.005194  &  0.037634 &  0.402861  \\
6 &            &            &            &            &  0.030379 &  0.383125  \\
7 &            &            &            &            &  0.030328 &  0.183904 \\
8 &            &            &            &            &           &  0.002054  \\
 \hline
 \end{tabular}
 \caption{Relative errors $E^{n,\eps}=|\frac{\sigma_{exact}- \sigma^{n,\eps}}{\sigma_{exact}}|$
 for iteration number $n=1,2,3,...$, and various $\eps$. The prescribed $\sigma$ is 5.  Empty space means the numerical result has already converged.}
 \label{ConvergenceOfDiffusionLengthSigmaCase1}
\end{table}

\begin{table}
\centering
 \begin{tabular}{|c |c |c |c |c |c |c| }
 \hline
 $n$ & $ \eps = 0.01 $ & $\eps = 0.02  $ &$\eps = 0.04  $ &$\eps = 0.08  $ & $\eps = 0.16  $ & $\eps = 0.32  $ \\
 \hline
1 &  0.677755 &   0.689595 &   0.737990 &   0.908029 &   0.414439   & 1.476995 \\
2 &  0.392867 &   0.408261 &   0.476915 &   0.160001 &   0.084646   & 0.691197  \\
3 &  0.089478 &   0.093146 &   0.092276 &   0.030613 &   0.029504   & 0.487247  \\
4 &  0.006387 &   0.007161 &   0.008500 &   0.008827 &   0.027198   & 0.158495  \\
5 &  0.000066 &   0.000377 &   0.002147 &   0.008453 &   0.027194   & 0.034154  \\
6 &  0.000033 &   0.000340 &   0.002115 &   0.008453 &              & 0.021585 \\
7 &           &            &            &            &              & 0.021471  \\
 \hline
 \end{tabular}
 \caption{Relative errors $E^{n,\eps}=|\frac{\sigma_{exact}- \sigma^{n,\eps}}{\sigma_{exact}}|$
 for iteration number $n=1,2,3,...$, and various $\eps$. The prescribed $\sigma$ is 10.}
 \label{ConvergenceOfDiffusionLengthSigmaCase2}
\end{table}

\begin{table}
\centering
 \begin{tabular}{|c |c |c |c |c |c |c| }
 \hline
 $n$ & $ \eps = 0.01 $ & $\eps = 0.02  $ &$\eps = 0.04  $ &$\eps = 0.08  $ & $\eps = 0.16  $ & $\eps = 0.32  $ \\
 \hline
1&   0.108867 &   0.109572 &   0.113283 &   0.126632 &   0.161406 &   0.237664 \\
2&   0.007695 &   0.008023 &   0.009952 &   0.016784 &   0.031044 &   0.040370  \\
3&   0.000070 &   0.000360 &   0.002080 &   0.008108 &   0.019861 &   0.024093  \\
4&   0.000031 &   0.000320 &   0.002038 &   0.008059 &   0.019782 &   0.023936 \\
5&            &            &            &            &   0.019782 &   0.023936 \\
 \hline
 \end{tabular}
 \caption{Relative errors $E^{n,\eps}=|\frac{\sigma_{exact}- \sigma^{n,\eps}}{\sigma_{exact}}|$
 for iteration number $n=1,2,3,...$, and various $\eps$. The prescribed $\sigma$ is 20.}
 \label{ConvergenceOfDiffusionLengthSigmaCase3}
\end{table}

\subsection{Validation of the diffusion-type model} Now, we are in the position to validate the diffusion model in estimating
the exciton diffusion length.  We are interested in identifying under which condition
the 1D model can be viewed as a good surrogate for the 2D model and how this condition relates to the property of organic semiconductors.

Again, only limited photoluminescence data from experiments are available and we have to solve
the forward model to generate data in our numerical tests. Specifically, given the exciton diffusion
length $\sigma$, the exciton generation function $G$, and the parametrization of the random interface
$h(\omega)$, we solve Eq. \eqref{eqn:rPDE1d} for a series of thicknesses $\{d_i\}$, and calculate the corresponding
expectations of the photoluminescence data $\{\tilde{\PL}_i\}$ according to Eq. \eqref{eqn:rPL1d}.
Therefore, $\{ d_i, \tilde{\PL}_i\}$ serves as the ``experimental'' data generated by the 1D model.
We then solve the minimization problem \eqref{min1} based on our numerically generated data $\{ d_i, \tilde{\PL}_i\}$ to estimate
the ``exact'' exciton diffusion length $\sigma$ in the presence of randomness, denoted by $\sigma_{exact}$ and will be used for
comparison.

\begin{figure}[h]
\centering
\subfigure[$\beta=0$]{\label{fig:CorrelationFunctionV2a}
\includegraphics[width=0.49\linewidth]{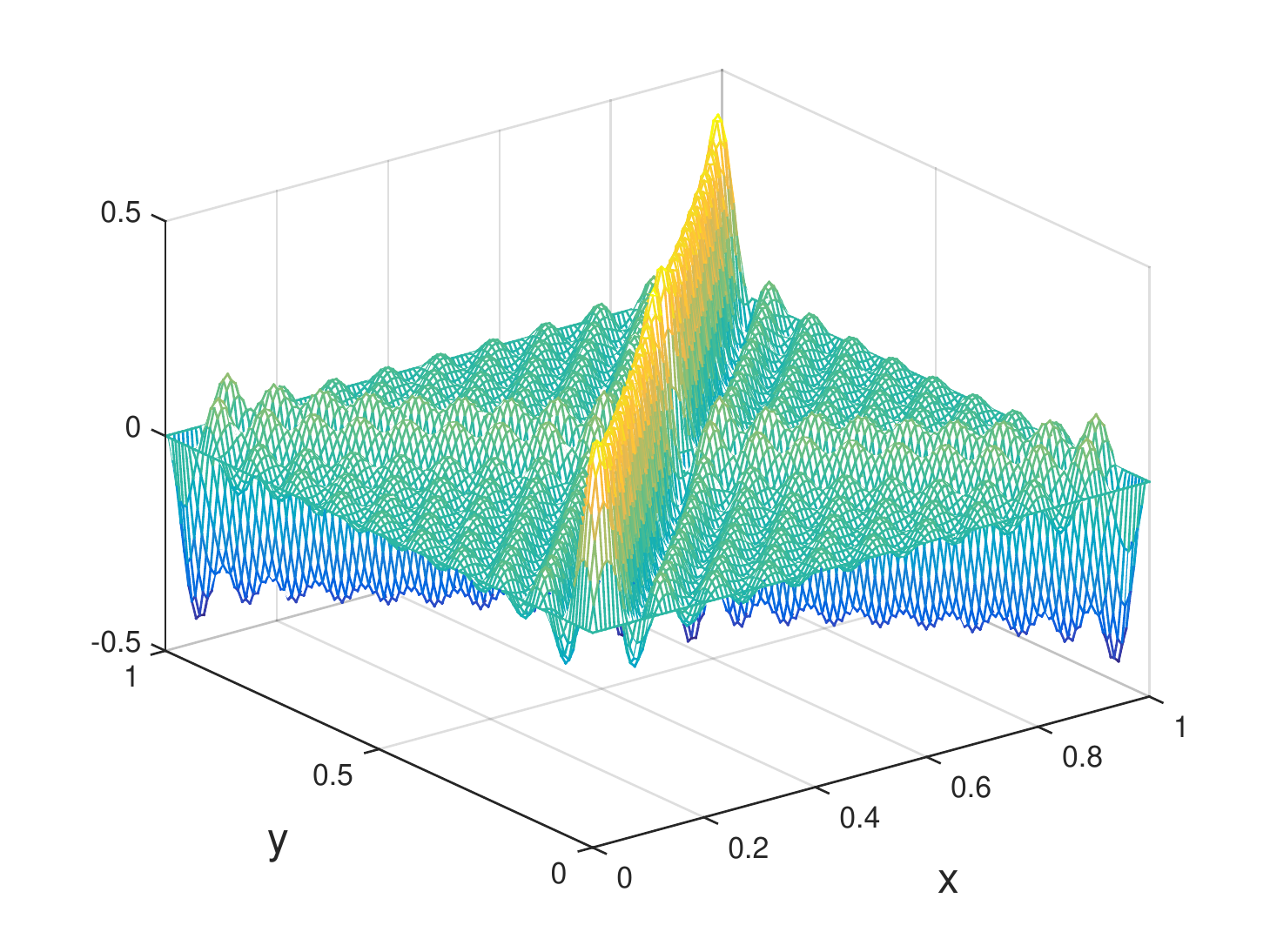}}%
\subfigure[$\beta=-2$]{\label{fig:CorrelationFunctionV2b}
\includegraphics[width=0.49\linewidth]{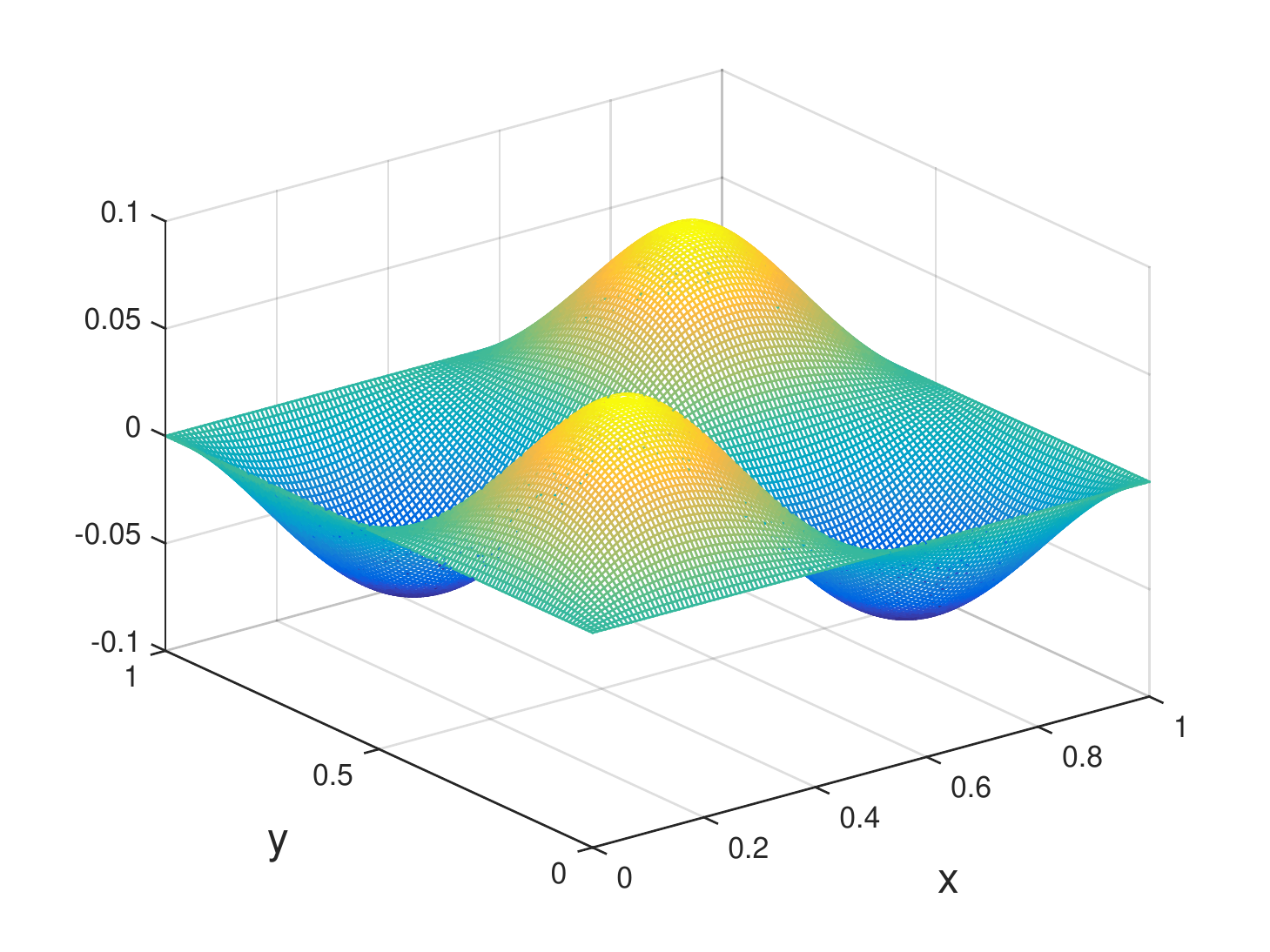}}%
\caption{\small  The covariance function of the random interface defined by Eq. \eqref{eqn:ex3_interface} for different $\beta$.
(a) $\beta=0$; (b) $\beta=-2$.}\label{fig:CorrelationFunctionV2}
\end{figure}

In our numerical tests, we use the 1D model \eqref{eqn:rPDE1d} with $\sigma=5$ and $\sigma=10$ to generate photoluminescence data.
$d_i=10i$, $i=1,...,10$, $\bar{h}=1$, and $\eps=\bar{h}/d_i$. We use $K=10$  random variables to parameterize the random interface.
We set $\lambda_{k}=k^{\beta}$, where $\beta\leq 0$ controls the decay rate of $\lambda_{k}$.
The random interface therefore takes the form
\begin{equation}\label{eqn:ex3_interface}
h(z,\omega) = \bar{h} \sum_{k=1}^{10}  k^\beta \theta_k(\omega)\sin(2k\pi\frac{z}{L})
\end{equation}
with $\theta_k(\omega)\sim U[-1,1]$. Figure \ref{fig:CorrelationFunctionV2} plots the covariance function of the
random interface defined by Eq. \eqref{eqn:ex3_interface} for $\beta=0$ and $\beta=-2$.
It is clear that the smaller the $\beta$, the larger the correlation length.

\begin{figure}[h]
\centering
\subfigure[$\sigma=5$]{\label{fig:Consis_EX3_U-11_DiffusionLengthSigma510Case1}
\includegraphics[width=0.49\linewidth]{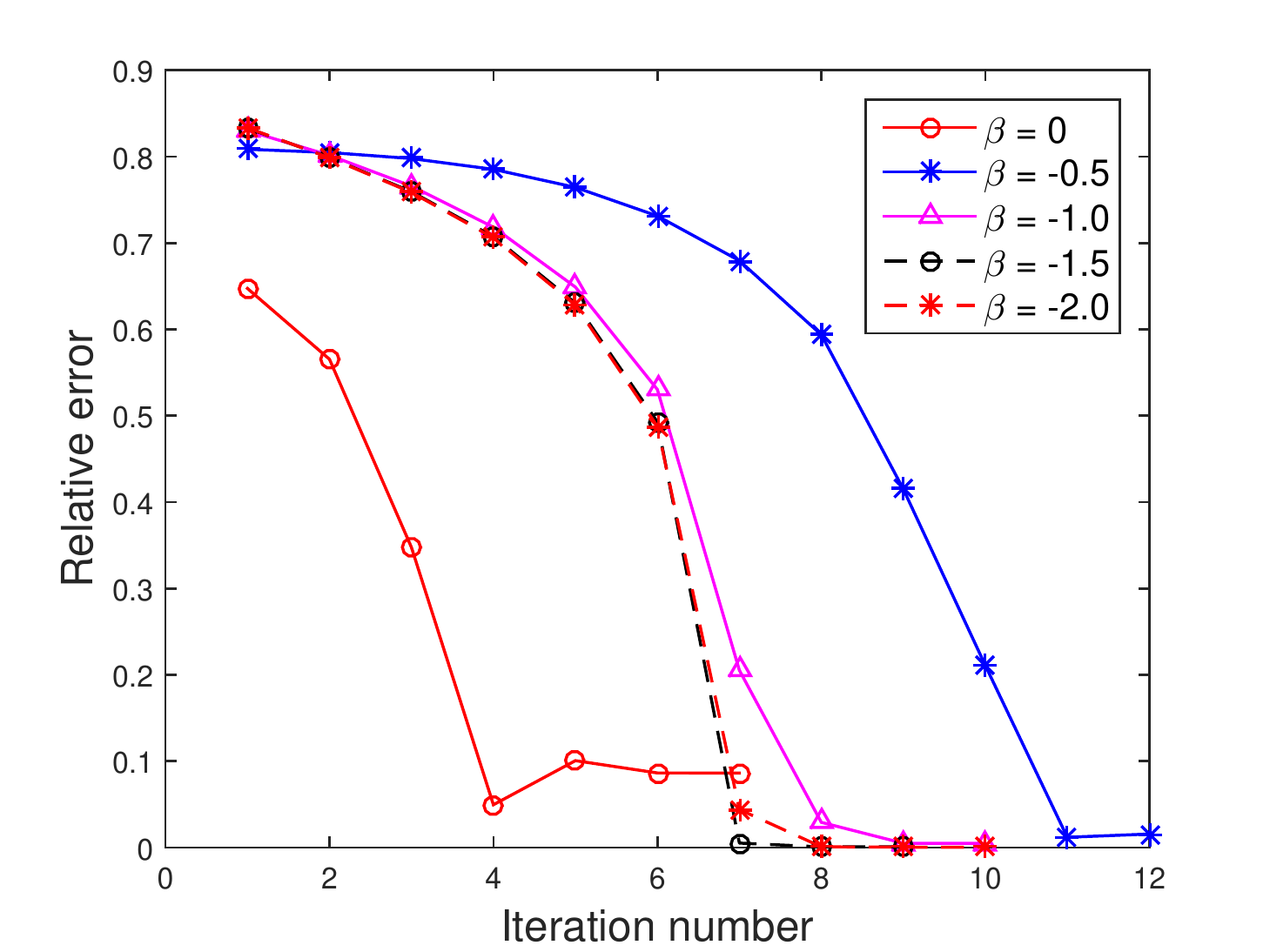}}%
\subfigure[$\sigma=10$]{\label{fig:Consis_EX3_U-11_DiffusionLengthSigma510Case2}
\includegraphics[width=0.49\linewidth]{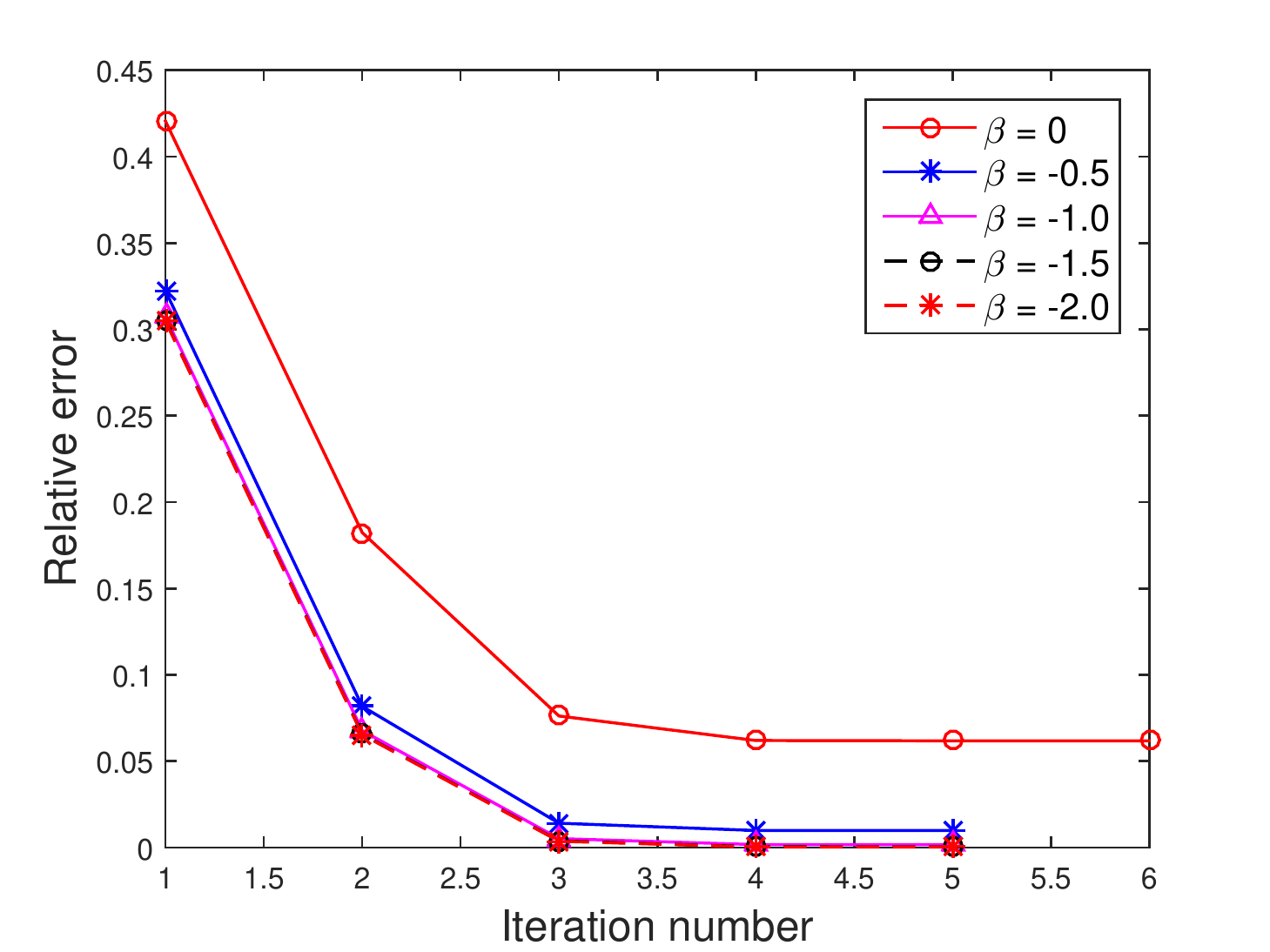}}%
\caption{\small Convergence history of the exciton diffusion length for various $\beta$, measured in
the relative error defined as $E^{n,\beta}=|\frac{\sigma_{exact}- \sigma^{n,\beta}}{\sigma_{exact}}|$
with $n$ the iteration number. The ``exact'' data is obtained by the 1D model (Eqs. \eqref{eqn:rPDE1d} and \eqref{eqn:rPL1d}) with
a prescribed $\sigma$. (a) $\sigma=5$; (b) $\sigma=10$.}\label{fig:Consis_EX3_U-11_DiffusionLengthSigma510}
\end{figure}

The convergence history of the exciton diffusion length for various $\beta$ is plotted in Figure \ref{fig:Consis_EX3_U-11_DiffusionLengthSigma510},
where the photoluminescence data is generated by the 1D model (Eqs. \eqref{eqn:rPDE1d} and \eqref{eqn:rPL1d}) with $\sigma=5$ and $\sigma=10$.
Again, the relative error is defined as $E^{n,\beta}=|\frac{\sigma_{exact}- \sigma^{n,\beta}}{\sigma_{exact}}|$, where $n$ is the iteration number, $\sigma_{exact}$ is the ``exact''  exciton diffusion length, and $\sigma^{n,\beta}$ is the numerical result defined in Eq. \eqref{estimate_sigma}. Note that $\sigma^{n,\beta}$ depends also on $\eps$ implicitly but we omit its dependence for convenience.
Tables \ref{Consis_EX3_SigmaCase5} and \ref{Consis_EX3_SigmaCase10} list the relative errors of our method for plotting Figure \ref{fig:Consis_EX3_U-11_DiffusionLengthSigma510}. The same criteria $|\ld^{(n)}-\ld^{(n-1)}|<10^{-4}$ is used here. The numerical exciton diffusion length obtained by our method converges to the reference one with the relative error less than $1\%$ when $\beta\leq -1$.

\begin{table}[h!]
\centering
 \begin{tabular}{|c |c |c |c |c |c |}
 \hline
 $n$ & $ \beta = 0 $ & $\beta = -0.5 $ &$\beta = -1.0$ &$\beta = -1.5 $ & $\beta = -2.0$  \\
 \hline
1 &   0.647340&   0.808503&   0.829974&   0.832905&   0.833375\\
2 &   0.565447&   0.804754&   0.801511&   0.798962&   0.798645\\
3 &   0.347626&   0.797750&   0.765892&   0.759374&   0.758543\\
4 &   0.049548&   0.785281&   0.718355&   0.707335&   0.705927\\
5 &   0.100595&   0.764404&   0.648943&   0.630438&   0.628049\\
6 &   0.086369&   0.731288&   0.530314&   0.492754&   0.487726\\
7 &   0.086302&   0.679544&   0.205261&   0.005156&   0.044032\\
8 &           &   0.593736&   0.029160&   0.000781&   0.000829\\
9 &           &   0.415006&   0.004952&   0.000777&   0.000410\\
10&           &   0.210863&   0.004758&           &   0.000410\\
11&           &   0.015892&           &           &           \\
12&           &   0.011628&           &           &           \\
 \hline
 \end{tabular}
 \caption{Relative errors $E^{n,\beta}=|\frac{\sigma_{exact}- \sigma^{n,\beta}}{\sigma_{exact}}|$
 for iteration number $n=1,2,3,...$, and various $\beta$. The prescribed $\sigma$ is 5.}
 \label{Consis_EX3_SigmaCase5}
\end{table}

\begin{table}[h!]
\centering
 \begin{tabular}{|c |c |c |c |c |c |}
 \hline
 $n$ & $ \beta = 0 $ & $\beta = -0.5 $ &$\beta = -1.0$ &$\beta = -1.5 $ & $\beta = -2.0$  \\
 \hline
1&   0.420520&   0.322303&   0.307239&   0.305058&   0.304686\\
2&   0.182350&   0.081737&   0.067669&   0.065660&   0.065316\\
3&   0.076150&   0.014086&   0.005175&   0.003874&   0.003647\\
4&   0.061968&   0.009871&   0.001699&   0.000493&   0.000281\\
5&   0.061776&   0.009857&   0.001690&   0.000483&   0.000272\\
6&   0.061776&           &           &           &           \\
 \hline
 \end{tabular}
 \caption{Relative errors $E^{n,\eps}=|\frac{\sigma_{exact}- \sigma^{n,\beta}}{\sigma_{exact}}|$
 for iteration number $n=1,2,3,...$, and various $\beta$. The prescribed $\sigma$ is 10.}
 \label{Consis_EX3_SigmaCase10}
\end{table}

Our numerical results show that a faster decay of the eigenvalues $\lambda_{k}$ leads to a better agreement between
the results of the 1D model and the 2D model. The smaller the $\beta$, the better the agreement. On the other hand,
the smaller the $\beta$, the larger the correlation length. Therefore, the larger the correlation length,
the better the agreement. Our study sheds some light on how to select
a model as simple as possible without loss of accuracy for describing exciton diffusion in organic materials.
In the chemistry community, it is known that under careful fabrication conditions \cite{DirksenRing:1991, Rodetal:2013},
organic semiconductors, including small molecules and polymers, can form crystal structures, which have
large correlation lengths. As a consequence, exciton diffusion
in these materials can be well described by the 1D model \cite{PetterssonRomanInganas:1999, Linetal:2013, Tamai:2015}. For organic materials with low crystalline order, i.e., small correlation length, however, our result suggests that the 1D model is not a good surrogate of the high dimensional models.

\section{Conclusion}
\label{sec:conclusion}

In this paper, we model the exciton diffusion by a diffusion-type equation with appropriate boundary conditions
over a random domain. The exciton diffusion length is extracted via minimizing the mean square error between
the experimental data and the model-generated data. Since the measurement uncertainty for the domain boundary
is much smaller compared to the device thickness, we propose an asymptotic-based method as the forward solver.
Its accuracy is justified both analytically and numerically and its efficiency is demonstrated by comparing
with the SC method as the forward solver. Moreover, we find that the correlation length
of randomness is the key parameter to determine whether a 1D surrogate is sufficient
for the forward modeling.

The discussion here focuses on the photoluminescence experiment. For the photocurrent experiment, from the modeling
perspective, the forward model is the same but the objective function is different. An exciton either contributes
to the photoluminescence or the photocurrent, so the photocurrent is defined as the difference between a
constant (total exciton contribution) and the photoluminescence \cite{Chen:2016}. Therefore, the proposed method can be applied straightforwardly with very little modification.

\medskip
\noindent \textbf{Acknowledgment.}
We thank Professor Carlos J. Garc\'{i}a-Cervia and Professor Thuc-Quyen Nguyen
for stimulating discussions. Part of the work was done when J. Chen
was visiting Department of Mathematics, City University of Hong Kong. J. Chen
would like to thank its hospitality. J. Chen acknowledges the financial support by
National Natural Science Foundation of China via grant 21602149. L. Lin and X. Zhou
acknowledge the financial support of Hong Kong GRF (109113, 11304314, 11304715).
Z. Zhang acknowledges the financial support of Hong Kong RGC grants (27300616, 17300817) and
National Natural Science Foundation of China via grant 11601457.
\section*{References}
\bibliographystyle{amsplain}

\providecommand{\MR}{\relax\ifhmode\unskip\space\fi MR }

\providecommand{\href}[2]{#2}


\newpage
\appendix

\section{Newton's method}
\label{sec:Newton}

The Newton's method works as follows: Given $\ld^{(0)}$, for $k=1,2,\ldots,$
\begin{equation}
\label{Newton}
\ld^{(k)} = \ld^{(k-1)} - \alpha_k \dfrac{\frac{\partial}{\partial \ld} J(\ld^{(k-1)})}
{\frac{\partial^2 }{\partial \ldsq}J(\ld^{(k-1)})},
\end{equation}
where $\alpha_k\in (0,1]$ is given by the line search technique.

For example, we take \eqref{min1} as the minimization problem and the domain mapping formulation
in \secref{sec:gpc} as the forward problem. Other combinations can be worked out similarly.
In 2D, for the first derivatives, we have
\begin{equation*}
\frac{\partial}{\partial \ld} J(\ld) = \frac2N \sum_{i=1}^N \left(\EE[\PL(\ld, d_i)] -
 \wt{\PL}_i\right)\EE[\frac{\partial \PL}{\partial \ld}]
\end{equation*}
and
\begin{equation}
\frac{\partial \PL(\ld,d_i)}{\partial \ld} = \int_0^1\int_0^1 \frac{\partial u}{\partial \ld}(d_i-h)\text{d}y\text{d}z.\nn
\end{equation}
Denote the derivatives of $u(y, z)$  with respective to the parameter $\ld$ by
\[
u_1(y, z) :=  \frac{\partial u}{\partial \ld} (y, z), ~~\mbox{ and }
u_2(y, z) :=  \frac{\partial^2 u}{\partial \ld^2} (y, z).
\]
Differentiating \eqref{eqn:PDE2d} with respect to $\ld$ directly,
we have
 \begin{equation}\label{eqn:u1PDE}
\ldsq\mathcal{L}u_1- u_1
= -2\ld\mathcal{L}u,   \quad (y,z)\in \dom_{\textrm{s}},
\end{equation}
and $u_1$ shares the same boundary condition as $u$.

For the second derivatives, we have
\begin{equation*}
\frac{\partial^2 }{\partial \ld} J(\ld) = \frac2N \sum_{i=1}^N \left(\EE[\frac{\partial \PL}{\partial \ld}]\right)^2 + \frac2N \sum_{i=1}^N \left(\EE[\PL(\ld,d_i)] -
 \wt{\PL}_i\right)\EE[\frac{\partial^2 \PL}{\partial \ld^2}]
\end{equation*}
and
\begin{equation*}
\frac{\partial^2 \PL(\ld,d_i)}{\partial \ldsq} = \int_0^1\int_0^1 u_2(y, z)(d_i-h)\text{d}y\text{d}z,\nn
\end{equation*}
and $u_2$ satisfies
\begin{equation}
\label{eqn:u2PDE}
\ldsq \mathcal{L} u_2 - u_2 = -2 \mathcal{L} u - 4\ld \mathcal{L}  u_1, \quad (y,z)\in \dom_{\textrm{s}}.
\end{equation}
Again, the same boundary condition applies for $u_2$.

To ease the implementation, we rewrite \eqref{eqn:u1PDE} and \eqref{eqn:u2PDE}
using \eqref{eqn:PDE2d}
\begin{equation}
\label{eqn:u12}
\begin{split}
\ldsq\mathcal{L}u_1- u_1 = -\frac{2}{\ld}(u-g), \\
\ldsq \mathcal{L} u_2 - u_2  = \frac{6}{\ldsq}(u-g) - \frac{4}{\ld}u_1.
\end{split}
\end{equation}

In the $k-$th step of Newton's method, knowing $\ld^{(k-1)}$, we solve
\eqref{eqn:PDE2d} and \eqref{eqn:bc2d} for $u^{(k-1)}$,
solve \eqref{eqn:u12} for $u_1^{(k-1)}$ and $u_2^{(k-1)}$,
and then update $\ld^{(k)}$ according to \eqref{Newton}.

In 1D, we have
\begin{equation*}
\begin{split}
\frac{\partial \PL(\sigma, d_i)}{\partial \ld} =  (d_i-\xi)\int_0^1 u_1(y)\text{d}y,\\
\frac{\partial^2 \PL(\sigma, d_i)}{\partial \ldsq} = (d_i-\xi)\int_0^1 u_2(y)\text{d}y
\end{split}
\end{equation*}
with $u_1(y)$ and $u_2(y)$ satisfying the same boundary condition as $u(y)$ (Eq. \eqref{eqn:bc1d})
and
\begin{equation*}
\begin{split}
\mathcal{L}_1 u_1(y) - u_1(y)  =  -\frac{2}{\ld}(u-G), \\
\mathcal{L}_1 u_2(y) - u_2(y)  =
\frac{6}{\ldsq} (u-G)-\frac{4}{\ld}u_1,
\end{split}
\end{equation*}
respectively.

\section{Asymptotic expansion}
\label{sec:derivation}

Using the change of variables, we first rewrite Eq. \eqref{eqn:rPDE2d} in $\tilde{x} = x/d$ and $\tilde{z} = z/L$
(still use $x$ and $z$ to represent $\tilde{x}$ and $\tilde{z}$). Note that the domain $\dom_\eps$ becomes
$\dom_{\textrm{s},\eps}:=\set{(x,z) \in (h(z,\omega)/d,1)\times(0,1)}$. Denote $\eps = \bar{h}/d$, then
\[
\dom_{\textrm{s},\eps}:=\set{(x,z) \in (\eps\tilde{h}(z,w),1)\times(0,1)},
\]
where $\tilde{h}(z,w) = \sum_k \lambda_k\theta_k(\omega)\phi_k(z)$.
Define $\mathcal{\tilde{L}} = \ldsq\left(d^{-2}\partial_{xx}+L^{-2}\partial_{zz}\right)-1$, then
\begin{numcases}
{ \label{eqn:rPDE2dc} }
 \mathcal{\tilde{L}} u_\eps(x,z) + g(x)   = 0,
&   $(x,z)\in \dom_{\textrm{s},\eps}$ \label{eqn:rPDE2d1c}
\\
\partial_x u_\eps(1, z) = 0, ~\ ~ u_\eps(\eps\tilde{h}(z,w), z)   = 0,  &   $ 0<z<1$ \label{eqn:rPDE2d2c}
\\
u_\eps(x, z)   =  u_\eps(x, z+1), &   $ \eps\tilde{h}(z,w) < x <1$ \label{eqn:rPDE2d3c}
\end{numcases}
with $u_\eps(x,z)$ and $g(x)$ representing $u(x,z)$ and $G(d-x)$ after the change of variables, respectively.

$\dom_{\textrm{s},\eps}$ depends on $\omega$, which brings great difficulty in numerical simulation.
It is easy to see, as $\eps\rightarrow 0$, $\dom_{\textrm{s},\eps}$ becomes a fixed domain $\dom_{\textrm{s}}=(0,1)\times(0,1)$.
To check the limit of $u_\eps$, we introduce  the following problem  for $u_\rex$ posed in the thin layer $L_\eps$:
\begin{equation}\label{eqn:ePDE2d}
\begin{cases}
& \tilde{\Lo} u_\rex+g(x)=0 \quad \text{in } L_\eps,\\
&u_\rex=u_\eps=0, \quad \partial_{\bn}  u_\rex=\partial_{\bn} u_\eps,\quad \text{on } \Gamma_\eps,\\
&u_\rex(x,z+L)=u_\rex(x,z), \quad \text{for }   (x,z) \in \overline{L_\eps}.
\end{cases}
\end{equation}
where $u_\eps$, the solution to equation \eqref{eqn:rPDE2dc}, is presumably given.
$\bn$ is the outward normal of $\dom_{\textrm{s}, \eps}$ on $\Gamma_\eps$.
At any point $(\eps \tilde{h}(z,\omega), z)\in \Gamma_\eps$, $\bn$ is parallel to the vector $(-1,\eps \tilde{h}'(z,\omega))$.
Here
\[
L_\eps=\{(x,z):0 \wedge (\eps \tilde{h}(z,\omega)) < x<  0 \vee \eps \tilde{h}(z,\omega),
0< z< 1\}
\]
and
\[
\Gamma_\eps := \overline{L_\eps}\cap\overline{\dom_{\textrm{s},\eps}}=\{(\eps \tilde{h}(z,\omega),z): 0\leqslant z\leqslant 1\}.
\]
In Figure \ref{fig:dom}, $L_\eps=\{(x,z):0 < x<  \eps \tilde{h}(z,\omega),
0< z< 1\}$ for positive $\tilde{h}(z,\omega)$ along $\Gamma_\eps$. For later use, we define
$\Gamma_0 :=\{(0,z): 0\leqslant z\leqslant 1\}$.

Note that Eq. \eqref{eqn:ePDE2d} is in fact a Cauchy problem of the time evolution equation
not a boundary-value problem of the elliptic PDE.
The velocity is specified on the interface $\Gamma_\eps$
by $\partial_{\bn}u_\eps$
and the wave travels along the normal $\bn$.
So, the solution of \eqref{eqn:ePDE2d} exists for $0\leq x\leq  \eps  \tilde{h}(z,\omega)$ \cite{Chen:2017}.
Particularly, we have  the existence of the value of $u_\rex$ at $x=0$.

\subsection{The solution on regular domain and its asymptotic expansion}
Now the solutions $u_\eps$ and $u_\rex$ are both well-defined on $\overline{\dom_{\textrm{s},\eps}}$ and $
\overline{L_\eps}$ by
\eqref{eqn:rPDE2dc} and \eqref{eqn:ePDE2d}, respectively.
In the next,  we introduce a function
 piecewisely defined by these two functions on the regular domain $\dom$ and want to find the correct equation
 for this function on $\dom$ in order to carry our asymptotic method.

Let $w_\eps$ be defined  on $\overline{\mathcal{D}}=\overline{\dom_{\textrm{s},\eps}} \cup \overline{L_\eps}$
as follows
 \begin{equation}
 w_\eps(x,z) := \begin{cases}
 u_\eps(x,z) &  \mbox{ in } \overline{\dom_{\textrm{s},\eps}},\\
 u_\rex(x,z) &  \mbox{ in }  \overline{L_\eps}.
 \end{cases}
 \end{equation}
This definition is justified by \eqref{eqn:ePDE2d} and  immediately implies the following obvious but important fact
which arises from the  boundary condition on the interface $\Gamma_\eps$ of $u_\eps$:
\begin{equation}\label{wonif}
w_\eps (x,z)=0 ~~~\mbox{ on } \Gamma_\eps.
\end{equation}

It is easy to see that $w_\eps$ is the unique solution to the following problem
where $u_\rex$ at $x=0$ is   given  {\it a prior} :
\begin{equation}\label{truesol}
\begin{cases}
&\tilde{\Lo}  w_\eps+G(d-x,z)=0 \quad \text{in } \mathcal{D},\\
&   w_\eps (0,z)= u_\rex(0,z),~~~\text{ for } 0\leqslant z\leqslant 1, \\
& \partial_x   w_\eps(1,z)=0, ~~~\text{ for } 0\leqslant z\leqslant 1,\\
&w_\eps(x,z+1)=w_\eps(x,z), \quad \text{for }   (x,z) \in \overline{\mathcal{D}}.
\end{cases}
\end{equation}

We start with  the  following ans\"{a}tz for $w_\eps$,
\begin{equation}\label{ans}
w_\eps(x,z)=\sum_{n=0}^\infty \eps^nw_{n}(x,z) \quad \text{for }(x,z)\in \overline{\dom}.
\end{equation}
Plug this ans\"{a}tz into  the equation \eqref{truesol}, and
match the terms at the same order of $\eps$,   then we obtain the following equations for   $w_n$
in $\dom$:
\begin{equation}
\begin{cases}
\tilde{\Lo} w_0+g(x)=0,      & \\
\tilde{\Lo} w_n=0,   ~~\quad n\geq 1. &
\end{cases}
\end{equation}
Next, we discuss the boundary conditions for these PDEs.
The two of the boundary conditions in \eqref{truesol},
 $\partial_x w_\eps(1,z)=0$ and $w_\eps(x,z+1)=w_\eps(x,z)$,
 do not depend on $u_\rex$. Thus, the ans\"{a}tz \eqref{ans} simply
  gives us the same boundary conditions for each $w_n$:
\begin{equation}
\partial_x w_n(1,z)=0, ~~~ \text{and }  w_n(x,z+1)=w_n(x,z).
\end{equation}
The  boundary condition of  \eqref{truesol}
at $x=0$, i.e.,  on $\Gamma_0$, depends on the data $u_\rex$ on this boundary.
If one works on this boundary condition,
it is possible to solve the Cauchy problem \eqref{eqn:ePDE2d}
for small $\eps$ analytically so that $u_\rex(x=0,z)$ can
be obtained in terms of $u_\eps$ (i.e., $w_\eps$),
and eventually certain connections for $w_n$ can be built.
But the use of the very original boundary condition \eqref{wonif} on $\Gamma_\eps\subset \partial \dom_{\textrm{s},\eps}$, not on $\partial \dom$, actually significantly simplifies the calculations
and  finally offers more friendly results.
The details follow below.

For the condition \eqref{wonif} on the interface $\Gamma_\eps$ where $x=\tilde{h}(z,\omega)$,
 \eqref{ans} implies
\begin{equation}\label{eqn:wn0}
w_\eps(\eps \tilde{h},z)=\sum_{n=0}^{\infty}\eps^n w_n(\eps \tilde{h},z)=0.
\end{equation}
The Taylor expansion in $\eps$
\begin{equation}\label{Taylor}
w_n(\eps \tilde{h},z)=\sum_{k=0}^\infty \frac{\eps^k\tilde{h}^k}{k!}\partial_x^k w_n(0,z),
\end{equation}
then gives \[
\sum_{k=0}^\infty \sum_{n=0}^{\infty}\eps^{n+k} \frac{\tilde{h}^k}{k!}\partial_x^kw_n(0,z)=0,
\]
which, by a change of the indices $m=k+n$, is equivalent to
\[
\sum_{m=0}^{\infty}\eps^{m} \sum_{k=0}^m\frac{\tilde{h}^k}{k!}\partial_x^kw_{m-k}(0,z)=0.
\]

Then by matching the terms with the same order of  $\eps$, we obtain:
\[
\sum_{k=0}^m \frac{\tilde{h}^k}{k!}\partial_x^kw_{m-k}(0,z)=0,\]
i.e.,
\begin{equation}\label{eqn:wm}
\begin{cases}
w_0(0,z) = 0, \\
w_{m}(0,z)=-\sum_{k=1}^m \frac{\tilde{h}^k}{k!}\partial_x^kw_{m-k}(0,z), ~~\forall m\geq 1.
\end{cases}
\end{equation}
This provides a recursive expression of  the boundary condition at $x=0$ for
the $m$-th order term $w_m$.

In summary,
the expansion of $u_\eps$ inside $\dom_{\textrm{s},\eps}$
is realized via  the expansion \eqref{ans},  $w_\eps = \sum_{n=0}^\infty w_n$, inside $\dom$.
Formally, each term $w_n$  satisfies the equation where the boundary condition   at $\Gamma_0\subset \partial \dom$
is defined recursively:
\begin{equation}\label{eqn:w0r}
\begin{cases}
&\tilde{\Lo} w_0+ g(x)=0 \quad \text{in } \dom,\\
& w_{0}(0,z)=0,  \quad \text{on }  \Gamma_0,
\\
&\partial_x w_0(1,z)=0, \quad  \text{for }  0\leqslant z\leqslant 1, \\
&w_0(x,z+1)=w_0(x,z), \quad  \text{for } (x,z) \in \overline{\dom},
\end{cases}
\end{equation}
and for $n\geq 1$,
\begin{equation}\label{eqn:wn}
\begin{cases}
&\tilde{\Lo} w_n=0 \quad \text{in } \dom,\\
& w_{n}(0,z)=-\sum_{k=1}^n \frac{\tilde{h}^k}{k!}\partial_x^kw_{n-k}(0,z),\quad \text{on }  \Gamma_0, \\
&\partial_x w_n(1,z)=0, \quad  \text{for }  0\leqslant z\leqslant 1, \\
&w_n(x,z+1)=w_n(x,z), \quad  \text{for } (x,z) \in \overline{\dom}.
\end{cases}
\end{equation}
In particular for  $m=1,2,3$, the above boundary conditions  on $\Gamma_0$
are
\begin{align}
& w_1(0,z)=- \tilde{h}\partial_x w_0(0,z),\\
& w_2(0,z)=-\tilde{h}\partial_x w_1(0,z)-\frac 12\tilde{h}^2\partial_{xx} w_0(0,z),
\label{w20z}\\
& w_3(0,z)=-\tilde{h}\partial_x w_2(0,z)-\frac 12\tilde{h}^2\partial_{xx} w_1(0,z)-\frac{1}{6} \tilde{h}^3\partial_{x}^3 w_0(0,z).
\end{align}

If we reverse the change of variables $\tilde{x} = x/d$ and $\tilde{z} = z/L$, \eqref{eqn:w0r} recovers \eqref{eqn:w0},
\eqref{eqn:wn} when $n=1$ and $n=2$ recovers the equations in \eqref{eqn:w1} and \eqref{eqn:w2}. Boundary conditions
in \eqref{eqn:w1} and \eqref{eqn:w2} can be recovered by using the inverse Lax-Wendroff procedure \cite{Chen:2017}.
In the boundary conditions for $w_n$ on $\Gamma_0$,  the
second and higher order partial derivatives with respect to $x$ may be converted to the partial derivatives with respect to $z$ by
repeatedly using the partial differential equations
\[
\sigma^2\partial_{xx}w_n+\sigma^2\partial_{zz}w_n-w_n
+\delta_{0,n}G(d-x)=0.
\]
Let us take order $n=0$ for example. Since $w_0(0,z)=0$, we have
\[
\sigma^2\partial_{xx}w_0(0,z)+\sigma^2\partial_{zz}w_0(0,z)-w_0(0,z)
+G(d)=\sigma^2\partial_{xx}w_0(0,z)+ G(d)=0,
\]
then
\[
\partial_{xx}w_0(0,z)=-\frac{1}{\sigma^2} G(d).
\]
This simplifies
\eqref{w20z}  to be
\begin{equation}
w_2(0,z)=-\tilde{h}(z)\partial_x w_1(0,z)+\frac{\tilde{h}^2(z)}{2\sigma^2} G(d).
\end{equation}
It is also easy to see that
\[
\partial_x^2 w_1(0,z) = - \partial^2_z w_1(0,z) + w_1(0,z)/\sigma^2.
\]
To compute $\partial_{x}^3 w_0(0,z)$, we take the derivative with respect to $x$ on both sides of the equation and get
\[
\sigma^2\partial_{x}^3w_0+\sigma^2\partial_{zz}\partial_xw_0-\partial_x w_0
-G'(d-x)=0,
\]
then taking values at $x=0$ yields
\[
\partial_{x}^3w_0(0,z)=\frac{1}{\sigma^2}\left[-\sigma^2\partial_{zz}\partial_xw_0(0,z)+\partial_x w_0(0,z)
+G'(d)\right].
\]
The other high order partial derivatives with respect to $x$ can also be converted to the partial derivatives with respect to $z$ similarly by using the corresponding partial differential equations.

\end{document}